\documentclass[11pt]{article}
\usepackage[utf8]{inputenc}
\usepackage{hyperref,amsmath,amssymb,natbib,enumerate,amsbsy,amsfonts, xcolor, graphics,graphicx, comment, algpseudocode, algorithm, bbm, booktabs, multirow,url, makecell, float, mathtools, pgfplots}
\RequirePackage{amsthm,amsmath,amsfonts,amssymb}
\usepackage[tableposition=above]{caption}
\usepackage{pifont}
\usepackage{tikz}

\usepackage{romannum}
\newtheorem{lemma}{Lemma}
\newtheorem{proposition}{Proposition}
\newtheorem{theorem}{Theorem}
\newtheorem{definition}{Definition}
\newtheorem{remark}{Remark}

\renewcommand{\theequation}{
\arabic{equation}
}
\makeatother

\usepackage{cleveref}

\pgfplotsset{compat=1.18}

\usepackage{xcolor}
\hypersetup{
    colorlinks,
    linkcolor={blue},
    citecolor={blue},
    urlcolor={black}
}
\begin{document}

\pagenumbering{arabic}

\title{Ask for More Than Bayes Optimal: A Theory of Indecisions for Selective Hypothesis Testing}
\author{Mohamed Ndaoud$^1$, Peter Radchenko$^2$, and Bradley Rava$^2$} 
\date{}

\footnotetext[1]{Department of Decision Sciences, ESSEC Business School, ndaoud@essec.edu}
\footnotetext[2]{Discipline of Business Analytics, University of Sydney Business School, peter.radchenko@sydney.edu.au, bradley.rava@sydney.edu.au}

\maketitle

\begin{abstract}
Selective classification is a powerful tool for automated decision-making in high-risk scenarios, allowing classifiers to act only when confident and abstain when uncertainty is high. Given a target accuracy, our goal is to minimize the number of indecisions, which are observations that we do not automate. For difficult problems, the target accuracy may be unattainable without abstaining from making a decision. By using indecisions, we can target a misclassification rate below the Bayes error rate, while minimizing overall indecision mass.

We provide a characterization of the optimal  risk in selective classification, establishing continuity and monotonicity properties that enable optimal indecision selection. We revisit selective inference via the Neyman-Pearson testing framework, where indecision enables control of Type~II error given a fixed Type~I error probability. For both classification and testing, we propose a calibration method, and analyze the excess risk of plug-in classifiers and the excess indecision mass produced by accuracy-based calibration. In the binary Gaussian mixture model, we identify an exponent-level phase transition, showing that minimal indecision can yield near-optimal accuracy even under poor class separation. Experiments on Gaussian mixtures and real datasets illustrate how indecision can improve selective accuracy.

\end{abstract}


\section{Introduction}
We address the problem of controlling a classifier's accuracy at any user-specified level through selective classification, regardless of the problem's inherent difficulty. Traditional classification frameworks are designed to approximate the Bayes optimal error rate as closely as possible. However, with the growing deployment of artificial intelligence (AI) systems in automated, high-stakes decision-making, it is not enough for a classifier to emulate or be asymptotically close to the Bayes optimal classification rule. Practitioners often require a prescribed level of reliability. When this level cannot be obtained, a classifier should be allowed to admit which instances it cannot automatically classify without violating this prescribed level, through the use of indecisions. 

This perspective is important because when the underlying problem is truly difficult, achieving control over the decisions made by an automated decision-making system may be impossible. When the number of potential classes is large or when the distributions of these classes are close enough, even the Bayes classifier may have error well above a practitioner's target level of reliability. The only way to always ensure reliable automated classification is through the use of indecisions, which are observations that are intentionally excluded from automated classification because their inherent difficulty prevents the algorithm from achieving the desired level of reliability. This naturally forms the basis of a selective classification paradigm, which allows us to effectively facilitate human-AI interaction: confident cases are automated, while difficult cases are deferred to expert human review. 

\begin{figure}[!t]
    \centering 
    \scalebox{0.9}{%
    \input{Figures/Density_Accuracy_Gamma_Est}}
    \caption{\label{accurary_gamma.plot}\small The best accuracy attained by the Bayes classifier and by selective classification as the distance between the underlying class distributions increases. The parameter $\Delta$ measures the separation between the high- and low-risk classes; larger values correspond to easier classification problems.}
\end{figure}

We illustrate this in the first plot of Figure~\ref{accurary_gamma.plot}, where the data for two classes are generated from Gaussian distributions with varying levels of separation between their centers, represented by $\Delta$. Our goal is to target a selective error rate of $1\%$ among automated decisions. When the level of separation is small, the Bayes classifier is unable to achieve our desired level of accuracy, as represented by the middle plot of Figure~\ref{accurary_gamma.plot}. This tells us that any standard classifier that is forced to make decisions on all instances can never hope to achieve our desired level of reliability. In contrast, selective classification, through the use of indecisions, can attain the target in this population illustration for all levels of separation. As the problem becomes easier, the Bayes error rate falls below our target level, and the need for indecisions disappears. This demonstrates where selective classification is the most useful: the regime in which our target level of performance falls (potentially far) below the Bayes error rate. The cost of achieving this level of reliability is shown in the rightmost plot of Figure~\ref{accurary_gamma.plot}, where the smallest indecision mass needed to achieve our desired level of reliability is displayed.

\begin{figure}[!t]
    \centering 
    \scalebox{1}{
    \input{Figures/intro_plot_method_compare}}
    \caption{\label{intro_method_comparison.plot}\small Difficult classification scenario with a nonlinear decision boundary. Left: One draw from the diamond data-generating mechanism at $\Delta'=3$, where larger $\Delta'$ gives sharper class separation around the nonlinear Bayes boundary. The middle and right panels compare the selective error rate and proportion of indecisions needed to achieve this error rate for methods targeting the nominal level $\alpha=0.10$. Several methods approach the desired error level, but they differ substantially in how many observations they leave undecided.} 
\end{figure}

Once we allow the use of indecisions, controlling the error rate across all automated decisions is only part of the problem. One can trivially always achieve a very small error rate by declaring nearly all observations as indecisions. However, this is largely uninformative and shifts a potentially larger than necessary burden on human experts. Our goal therefore should be to achieve any user-specified error rate, while minimizing the number of observations identified as indecisions. 

Achieving this goal in practice is challenging. Some selective classifiers can attain a user-specified level of reliability but struggle to maintain it in the most difficult classification settings, most notably near a phase transition, where large improvements in accuracy can be realized with negligible amounts of indecisions. In contrast, other methods may fail to achieve the desired reliability altogether when the classification problem is sufficiently difficult or the true decision boundary is complex. Figure~\ref{intro_method_comparison.plot} illustrates this point in a nonlinear setting, where the two classes are separated by a diamond-shaped decision boundary. In this case, the difficulty of the problem is controlled by $\Delta'$, which plays a similar role to $\Delta$ in Figure~\ref{accurary_gamma.plot}. As $\Delta'$ increases, the posterior probability for each class changes more abruptly, causing the two classes to be more separated. 

The simulation uses $10{,}000$ observations split evenly between training and calibration samples, together with an additional $10{,}000$ observations for testing. For illustration, we compare two leading selective classification approaches adapted to target the overall error rate. The first, classification with reject using a hinge loss (orange / square), is tuned on an external calibration set with cross-validation to target a $10\%$ error rate \citep{bartlett08a}. The second is a conformal inference approach based on conformal $p$-values (blue / triangle) that classifies each observation independently into each class \citep{Jin_23}. Because this conformal procedure is adapted from its original objective, its original guarantee need not apply to the selective error considered here. Two features are immediately apparent. Classification with reject struggles to maintain the desired level of reliability and also uses substantially more indecisions than the oracle benchmark. The conformal $p$-value approach approaches the target selective error in this experiment, but it struggles near the phase transition of interest, which is of high value for practitioners. Moreover, both approaches are limited to the conditional error rate, and it remains unclear how they can be extended to accommodate more complex error metrics, such as Type~I and Type~II errors.\footnote{Classification with reject uses an SVM base classifier, whereas the conformal $p$-value method and our proposed approach use XGBoost.}

To address all of these concerns, this work aims to answer the following central question: how can one control the error rate of a classifier at any user-specified level while using the smallest possible number of indecisions? Our proposed approach (green / circle), in contrast to the competing methods, closely tracks the desired error rate in this experiment while nearly matching the minimum necessary proportion of indecisions.

The rest of this paper is structured as follows. In Sections~\ref{sec:problem_formulation} and~\ref{sec:lit_review_intro}, we formulate our problem and review related work. In Sections~\ref{sec:general} and~\ref{sec:gen_hyp_test}, we derive optimal selective classifiers for general binary classification and for a hypothesis-testing formulation involving Type~I and Type~II errors. We consider both a fixed indecision mass and fixed target error levels, and we provide calibration algorithms.

Section~\ref{sec:theory} studies the phase transition in the Gaussian mixture model and then analyzes plug-in rules based on an estimate of the conditional probability function. Section~\ref{appendix:sec:ext} develops extensions under the monotone likelihood-ratio property and for multiclass classification. Section~\ref{sec:simu} presents simulation studies of the phase transition and the calibration procedures. Building upon the simulation analysis, Section~\ref{sim:real_data_compas} applies the classification and hypothesis-testing procedures to COMPAS recidivism data \citep{Angwin_16}. Proofs of all theorems are in the appendix. 

\subsection{Problem Formulation}
\label{sec:problem_formulation}
We observe a random variable $X$ on a measurable space $(\mathcal{X},\mathcal{U})$ such that $X$ is distributed according to a mixture model, where with probability $p_1$ its probability measure is given by $P_{1}$ and with probability $p_{2}$ its probability measure is $P_{2}$. We assume $p_1,p_2>0$ and $p_1+p_2=1$; the degenerate case $P_1=P_2$ is allowed. Let $f_{1}$ and $f_{2}$ be densities of $P_{1}$ and $P_{2}$ with respect to some dominating measure that we will further denote by $\mu$. Denote by $Y$ the labeling quantity such that $Y=1$ if the distribution of $X$ is $P_{1}$ and $Y=2$ if it is $P_{2}$. We are interested in the problem of predicting the true label $Y$ with an estimator $\widehat{Y}$, with the quality of our estimate measured either conditional on making a decision through supervised classification or conditional on the true label $Y$ through controlling Type~I and Type~II errors in a hypothesis-testing formulation. 

As estimators of $Y$, we consider any measurable functions $\widehat Y=\widehat Y(X)$ taking values in $\{1,2\}$. 
Such estimators will be called {\it classifiers}. 
We define the loss of a classifier $\widehat Y$ as the indicator of whether a mistake is made, that is $\mathbf{1}\{ \widehat Y \neq Y\}$, where $\mathbf{1}(\cdot)$ is the indicator function. The performance  of~$\widehat Y$ is measured by its expected risk ${\mathbf P}_Y (\widehat Y \neq Y)$, also known as classification error or misclassification rate, or by ${\mathbf P}_Y (\widehat Y = Y)$, referred to as accuracy. 

We denote by ${\mathbf E}_Y$ the expectation with respect to probability measure ${\mathbf P}_Y$ of $X$ with labeling~$Y$. Observe that ${\mathbf P}_Y (\widehat Y \neq Y) = p_1 {\mathbf P}_1 (\widehat Y = 2) + p_{2} {\mathbf P}_{2} (\widehat Y = 1)$, which is a weighted sum of the Type~I and Type~II errors. The classical theory of classification gives a precise characterization of the Bayes risk, $\inf_{\tilde Y} {\mathbf P}_Y(\tilde Y \neq Y)$, where $\inf_{\tilde Y}$ denotes the infimum over all measurable classifiers. In particular, an optimal Bayes classifier is defined, with ties assigned to class~1, by
\[
Y^*(X)=
\begin{cases}
1, & p_1f_1(X)\ge p_2f_2(X),\\
2, & p_1f_1(X)<p_2f_2(X).
\end{cases}
\]
Moreover, the corresponding risk is given by
$$
\inf_{\tilde Y}  {\mathbf P}_Y(\tilde Y \neq Y) = {\mathbf P}_Y( Y^* \neq Y) = \int (p_1 f_1 \wedge p_{2} f_2) d\mu = \frac{1}{2} - \frac{1}{2}\int |p_1 f_1 - p_{2} f_{2}| d\mu.
$$
In particular, the Bayes risk is determined by the overlap between the two distributions. When $f_1$ and $f_2$ are close, any classifier performs poorly, which serves as motivation for the present work. Our goal is to introduce and study a framework in which accuracy arbitrarily close to one can be achieved with the help of indecisions, when the required posterior-confidence regions have positive probability.

In order to break the statistical barrier given by the Bayes risk, we allow our estimator a degree of freedom where it only makes a decision when it is sufficiently confident. Depending on the targeted accuracy level, the classifier may have to discard some of the observations. More precisely, given an \emph{indecision level} $\gamma$, we will consider the new risk:
\begin{equation}\label{eq:minimax risk}
  \mathcal{R}(\gamma) := \inf_{\tilde{Y}_{\gamma}}  {\mathbf P}_Y \left(\tilde{Y}_{\gamma} \neq Y | \tilde{Y}_{\gamma} \neq 0 \right),  
\end{equation}
where $\inf_{\tilde Y_\gamma}$ is the infimum over all classifiers taking values in $\{0,1,2\}$ such that ${\mathbf P}(\tilde Y_\gamma = 0) = \gamma$. In other words, we are interested in the lowest conditional misclassification risk when decisions are made for a pre-specified proportion of observations.

\subsection{Related Literature}\label{sec:lit_review_intro}
The concept of binary classification with indecisions has been well studied by different communities. It is known by several names, such as ``classification with a reject option'', ``selective classification'', ``no-decision classification'', ``classification with abstention'', and ``human--AI collaboration''. The corresponding approaches involve classifiers that are allowed not to make a decision when the class probabilities used for making a decision are too close to each other. For clarity, throughout this paper we refer to all of these areas of the literature under the umbrella term of selective classification \citep{Ran_10}.

Selective classification has traditionally encompassed two primary forms of observation rejection, referred to in this paper as indecision: ambiguity rejection and novelty rejection \citep{survey}. Ambiguity rejection arises when a model cannot confidently differentiate between two or more classes for a given observation \citep{Chow_57, Hellman_70, Fukunaga_72}. In contrast, novelty rejection applies to observations that cannot be reliably assigned to any predefined class \citep{Cordella_95, Seo_00, Vailaya_00}. 
In this work, we introduce a distinct rejection paradigm that shares characteristics with both ambiguity and novelty rejection but is fundamentally driven by classification accuracy. Specifically, we propose accuracy rejection, a paradigm that rejects observations that cannot be classified without the accuracy falling below a predefined threshold, irrespective of whether they exhibit ambiguity or novelty. This type of rejection has been explored in prior work under a variety of assumptions and guarantees \citep{Shekhar_19, Rava_25}. While accuracy rejection bears similarities to ambiguity rejection, it explicitly prioritizes the maintenance of a specified classification accuracy level.

\subsubsection{Classification with Reject}\label{sec:lit_review_cwr}

In the binary classification setting, the classification with reject paradigm seeks to incorporate indecisions into a classifier by optimally picking the cost of indecisions $d$ and the threshold $\delta$, then minimizing the modified cost function $\mathbf P\{H \cdot f(x) < -\delta\} + d \cdot \, \mathbf P\{|H \cdot f(x)| \leq \delta\},$ where $H:= 2 Y - 3$.
Although closely related to our proposed selective classification framework, classification with reject differs in several important aspects. In particular, existing reject-based approaches do not provide clear guidance about how to obtain rigorous accuracy guarantees while simultaneously minimizing the number of identified indecisions.

There have also been recent investigations of how best to incorporate selective classification into modern machine-learning algorithms through the lens of convex optimization \citep{Yuan_10}. The works of \cite{Grandvalet_08} and \cite{WEGKAMP_11} studied selective classification alongside support vector machines, while \cite{Cortes_16} investigated simultaneously learning a classifier and a selection rule. Recently, selective classification has been studied as a way to manage limited resources \citep{Valade_24}. The impact of plug-in classifiers on oracle selection rules has also been studied by \cite{Denis_20} and \cite{Lei_14} under continuity assumptions near the decision threshold. This work was further explored by \cite{Shekhar_19}, who controlled the abstention constraint with high probability while also addressing discontinuities in the empirical CDF. 

\subsubsection{Multiple Testing, Outlier Detection, and Conformal Inference}\label{sec:lit_review_calib}

Another stream of work looks at identifying a calibrated selected set of observations that focuses on controlling coverage over a smaller set of observations, up to a user-specified level \cite{Lei_14}. There are overlaps with conformal inference where the goal is to create prediction sets that contain the true classification label, up to a user-specified level of coverage \cite{Vovk_99, Vovk_05}. Other works have aimed to bridge the gap between overall coverage and calibrated decision-making. For binary classification, \cite{Lei_14} constructed confidence sets that can be calibrated for each class at a user-specified level. Another stream of literature offers calibrated decision-making through control over the False Selection Rate (FSR), the expected proportion of erroneous decisions among selected observations \citep{Gang_2022, Huo_24, Marandon_24, Jin_23, Rava_25, FASI_multi}. In a similar spirit, \cite{Sun_15} developed a decision-theoretic framework that utilized indecisions to control the FSR and \cite{Wang_24} recently used indecisions in the sequential setting.

On the modern application side, selective classification has also been used to address societal issues, such as fairness in decision-making. The works of \cite{Schreuder_21} and \cite{Rava_25} have independently investigated how to transform off-the-shelf classifiers into fair selective classifiers through empirical risk minimization and calibrated selection rules, respectively.

\subsection{Comparison to Existing Approaches}\label{sec:comparison_existing_approaches}

Our method is most similar to classification with reject and selection using conformal $p$-values. Classification with reject introduces an indecision option directly into the learning objective, while selection with conformal $p$-values aims at controlling the rate of erroneous decisions on a selected group of interesting observations. However, our method is motivated by constrained Bayes optimization rather than by coverage or rejection-cost objectives. This perspective enables theoretical guarantees in specific settings and leads to stronger performance in the most difficult regimes, such as the phase-transition setting studied later in this paper.
 
There are many conformal methods tailored for classification that seek to obtain marginal coverage for prediction sets rather than selective error control \citep{Romano_20, Angelopoulos_25, Sesia_25}. However, obtaining marginal coverage is a distinctly different objective than selective classification, as noted in \cite{Jin_23, Rava_25}. Providing prediction sets that contain the true label on average does not guarantee that the prediction sets that contain only one label are accurate. This distinction matters because selective classification evaluates error only on the data-dependent subset of observations assigned a definitive label.

\subsection{Our contributions} 
\label{sec:contributions}

We start with a characterization of the optimal selective risk \eqref{eq:minimax risk} in binary classification (Section~\ref{sec:general}), which we later generalize to multiclass classification (Section~\ref{sec:multi}). Our threshold characterization covers continuous score distributions and score plateaus whenever the required boundary mass is attainable. Along the way, we show that the map $\gamma \mapsto \mathcal{R}(\gamma)$ is continuous and non-increasing. In other words, for any given reachable level of accuracy, we can find a matching optimal indecision level and the corresponding classifier. These findings are extended to a hypothesis-testing problem in which, given a Type~I error, we wish to control the Type~II error. Our focus is the conditional-error, minimal-indecision formulation of this problem.

Sections \ref{sec:binary_calib} and \ref{sec:np_calib} are dedicated to our fully adaptive methodology. We offer novel analysis of both the classification and the hypothesis testing settings, with corresponding simulation and real data analyses in sections~\ref{sim:np_binary} and~\ref{sim:real_data_compas}. We explain how to calibrate the indecision region given a plug-in rule $\hat{\eta}$, in order to either achieve a level of accuracy or match a level of indecisions for both problems of classification and testing. The setting with a fixed budget of indecisions, which arises when practitioners can allocate only a limited number of observations for human review, remains relatively underexplored in selective inference and distinguishes our work from prior approaches.

In Section \ref{sec:gaussian}, we focus on the binary Gaussian mixture model given a fixed separation between the centers. We fully characterize the ``sharp'' phase transition of classification in terms of indecisions. It is well established that, in order to achieve a level of accuracy of order $1-\delta$, the separation between centers $\Delta$ has to be of the order $\sqrt{2\log(1/\delta)}$ where the constant $2$ is sharp. When the separation is of order $c\sqrt{2\log(1/\delta)}$ for some $c<1$, we need indecisions to reach the level of accuracy $1-\delta$. We give a sharp characterization of indecisions in this case. Interestingly, as long as $1/2<c<1$, we show that the optimal amount of indecisions is of order $o(1)$, meaning that by allowing only a negligible proportion of indecisions we can reach the level of misclassification~$\delta$ even in the case where the class distributions are not well-separated.  These findings are illustrated by numerical experiments in Section \ref{sim:phase_transition}. More generally, the optimal procedure is based on thresholding the likelihood ratio between distributions $f_1$ and $f_2$, which can be encoded through the regression function $\eta$. In practice, we can use a training sample to learn $\eta$. In Section \ref{sec.plugin}, we quantify the loss induced by the estimation of~$\eta$. First, under reasonable assumptions similar to the usual margin condition, we show that for a fixed level of indecisions, the accuracy of the plug-in procedure is comparable to that of the oracle and, in general, we can expect consistency of the plug-in approach. Second, we also show that if calibration is done with respect to the accuracy, i.e., if we tune the plug-in classifier to reach a given accuracy level, then the amount of indecisions is also controlled as the sample size grows, although not necessarily consistently.  

Finally, we offer two extensions of our theory in Section~\ref{appendix:sec:ext}. First, we emphasize the special case where the likelihood ratio $f_2/f_1$ satisfies the monotone likelihood-ratio property. In this setting, we do not need a training sample, as we can simply threshold the observations themselves instead of the scores $\eta$. This is typically the case for location models under log-concave distributions. We also show how to calibrate our procedure in this setting. Second, we extend our classification theory to the multiclass setting.

\subsection{Notation}

Throughout the paper we use the following notation. 
For given quantities $a_{n}$ and $b_n$, we write $a_{n}  \lesssim b_{n}$ ($a_{n} \gtrsim b_n$) when $a_{n} \leq c b_{n}$ ($a_{n} \geq c b_{n}$) for some absolute constant $c>0$. In the case $a_n/b_n \to 0$, we use the notation $a_n = o(b_n)$.  We also write $a_{n} \approx b_{n} $ if $a_{n} \lesssim b_{n}$  and $a_{n} \gtrsim b_n$. For any $a,b \in \mathbf{R}$, we denote by $a\vee b$ ($a \wedge b$) the maximum (the minimum) of $a$ and $b$.  Finally $c_0$, $c_1$, $c$ are used for positive constants whose values may vary from theorem to theorem.
We use $\alpha$ for an overall selective-risk target and $\alpha_1,\alpha_2$ for the Type~I and Type~II targets, respectively. A fixed proportion of indecisions is denoted by $\gamma$, and the smallest proportion of indecisions that meets the specified error target is denoted by $\gamma^*$.

\section{General Binary Classification}\label{sec:general}

In this setting, we want to find the \emph{smallest} possible indecision region that is able to control the accuracy (or, similarly, misclassification risk) of our classifier. By using the minimum necessary amount of indecisions, we can automate as many decisions as possible and send only the necessary exceptions for costly human review.

For a given classifier $\tilde{Y}$, our objective is to control the conditional misclassification risk, given that a decision has been made. We define the optimal selective risk as
$$
  \mathcal{R}(\gamma) := \inf_{\tilde{Y}_{\gamma}}  {\mathbf P}_Y \left(\tilde{Y}_{\gamma} \neq Y | \tilde{Y}_{\gamma} \neq 0 \right),  
$$
where the infimum is over all classifiers taking values in $\{0,1,2\}$ such that ${\mathbf P}(\tilde Y_\gamma = 0) = \gamma$. Our goal at the end of this section is to ensure that this conditional risk does not exceed a pre-specified threshold $\alpha$, while minimizing the number of indecisions. 
To accomplish this, we first study the comparatively easier setting where, for a fixed proportion of indecisions $\gamma$, we aim to find the best achievable level of accuracy. Understanding the case with a fixed amount of indecisions provides a framework for developing a procedure that targets a desired selective risk.

\subsection{Fixed proportion of indecisions}\label{sec:binary_class_fixed_indec}

Here, we focus on the binary case where we only have two classes. For a given level of indecisions $\gamma\in[0,1)$, we define the optimal indecision region~$\Theta_{\gamma}$, satisfying $\mathbf{P}_Y(\Theta_{\gamma}) = {\gamma}$. We will show that there exists a value $\tau_{\gamma} \in [1/2,1] $ such that
\begin{equation}\label{ind_reg.equn}
\Theta_{\gamma}:= \left\{ 1 - \tau_{\gamma} <  \frac{p_1f_1(X)}{p_1f_1(X) + p_{2}f_{2}(X)} < \tau_{\gamma} \right\} \cup \mathcal{M}_{\gamma},
\end{equation}
where $\mathcal{M}_{\gamma}$ is any subset of $\left\{   \frac{p_1f_1(X) \vee  p_2f_2(X)}{p_1f_1(X) + p_{2}f_{2}(X)} =  \tau_{\gamma} \right\} $ such that $\mathbf{P}_Y(\Theta_{\gamma}) = {\gamma}$.

We define $\eta(\cdot)$ as the conditional probability function $\eta(x) = \mathbf{P}(Y=1|X=x)$, that is,
\begin{equation}\label{cond_dens.equn}
    \eta(X) = \frac{p_1f_1(X)}{p_1f_1(X) + p_{2}f_{2}(X)}.    
\end{equation}

It is natural to observe that the optimal indecision region concentrates around where $\eta(X)$ is close to $1/2$. We note that our~threshold $\tau_{\gamma}$ plays a similar role to the constant~$d$ in \cite{Herbei_Wegkamp_06}. We also note that when $\eta(X) \vee (1-\eta(X)) = \tau_{\gamma}$, i.e., we are at the frontier of making an indecision, then we might randomly choose to reject or not. The Bayes oracle classifier with a ${\gamma}$ proportion of indecisions is given by
\begin{equation}\label{eq:oracle}
    Y_{\gamma}^*(X)=
    \begin{cases}
    0, & X\in\Theta_{\gamma},\\
    1, & X\notin\Theta_{\gamma}\ \text{and}\ p_1f_1(X)\ge p_2f_2(X),\\
    2, & X\notin\Theta_{\gamma}\ \text{and}\ p_1f_1(X)<p_2f_2(X).
    \end{cases}
\end{equation}
as shown in the following result.
\begin{theorem}\label{thm:acc:lower_bound}
   Given $\gamma$, the classifier $Y_{\gamma}^*$ is optimal for the risk $\mathcal{R}(\gamma)$. Moreover, we have that
   $$ 
   \mathcal{R}(\gamma) = {\mathbf P}_Y( Y_{\gamma}^* \neq Y |  Y_{\gamma}^* \neq 0 ) = \frac{\int_{\Theta_{\gamma}^c} (p_1f_1 \wedge p_{2}f_{2}) d\mu }{1-{\gamma}},
   $$
   where $\Theta_{\gamma}^c$ denotes the complement of the set $\Theta_{\gamma}$.
\end{theorem}
It follows from our proof that~$\tau_{\gamma}$ is an increasing function of~$\gamma$, and as $\gamma \to 0$, we recover the classical classification result without indecisions. 

On the one hand, when the random variable $\eta(X)$ has atoms and, in particular, $\mathbf{P}(\eta(X) \vee (1-\eta(X))  = \tau_{\gamma}) \neq 0$, the set $\mathcal{M}$ is nonempty, and we call the region $\mathcal{M}$ a ``plateau'', where the indecisions are picked randomly, as shown in the left panel of Figure~\ref{fig:combined_intro_plots}. On the other hand, if $\eta(X)$ has no atoms, then the indecision region is unique up to $\mu$-null sets.

We emphasize that the indecision region is not necessarily an interval, as illustrated in the rightmost plot of Figure~\ref{fig:combined_intro_plots}. Consequently, constructing the indecision region requires prior knowledge of the conditional probability~$\eta$.

 \begin{figure}[!t]
    \centering
    \begin{minipage}[t]{0.48\textwidth}
        \centering
        \scalebox{0.65}{%
\begin{tikzpicture}[x=1pt,y=1pt]
\definecolor{fillColor}{RGB}{255,255,255}
\path[use as bounding box,fill=fillColor,fill opacity=0.00] (0,0) rectangle (361.35,180.67);
\begin{scope}
\path[clip] (  0.00,  0.00) rectangle (361.35,180.67);
\definecolor{drawColor}{RGB}{255,255,255}
\definecolor{fillColor}{RGB}{255,255,255}

\path[draw=drawColor,line width= 0.5pt,line join=round,line cap=round,fill=fillColor] (  0.00,  0.00) rectangle (361.35,180.68);
\end{scope}
\begin{scope}
\path[clip] ( 42.21, 41.95) rectangle (356.35,175.67);
\definecolor{fillColor}{RGB}{255,255,255}

\path[fill=fillColor] ( 42.21, 41.95) rectangle (356.35,175.67);
\definecolor{fillColor}{RGB}{86,180,233}

\path[fill=fillColor] ( 59.25, 45.62) --
	( 64.25, 45.62) --
	( 64.25, 50.61) --
	( 59.25, 50.61) --
	cycle;

\path[fill=fillColor] ( 65.81, 45.73) --
	( 70.80, 45.73) --
	( 70.80, 50.73) --
	( 65.81, 50.73) --
	cycle;

\path[fill=fillColor] ( 72.36, 45.97) --
	( 77.35, 45.97) --
	( 77.35, 50.96) --
	( 72.36, 50.96) --
	cycle;

\path[fill=fillColor] ( 78.91, 46.45) --
	( 83.91, 46.45) --
	( 83.91, 51.45) --
	( 78.91, 51.45) --
	cycle;

\path[fill=fillColor] ( 85.46, 47.38) --
	( 90.46, 47.38) --
	( 90.46, 52.37) --
	( 85.46, 52.37) --
	cycle;

\path[fill=fillColor] ( 92.02, 48.99) --
	( 97.01, 48.99) --
	( 97.01, 53.98) --
	( 92.02, 53.98) --
	cycle;

\path[fill=fillColor] ( 98.57, 51.76) --
	(103.56, 51.76) --
	(103.56, 56.76) --
	( 98.57, 56.76) --
	cycle;

\path[fill=fillColor] (105.12, 56.02) --
	(110.12, 56.02) --
	(110.12, 61.01) --
	(105.12, 61.01) --
	cycle;

\path[fill=fillColor] (111.67, 62.50) --
	(116.67, 62.50) --
	(116.67, 67.50) --
	(111.67, 67.50) --
	cycle;

\path[fill=fillColor] (118.23, 71.49) --
	(123.22, 71.49) --
	(123.22, 76.48) --
	(118.23, 76.48) --
	cycle;

\path[fill=fillColor] (124.78, 83.26) --
	(129.77, 83.26) --
	(129.77, 88.26) --
	(124.78, 88.26) --
	cycle;
\definecolor{drawColor}{RGB}{213,94,0}

\path[draw=drawColor,line width= 0.4pt,line join=round,line cap=round] (133.83,100.04) circle (  2.50);

\path[draw=drawColor,line width= 0.4pt,line join=round,line cap=round] (140.38,116.46) circle (  2.50);

\path[draw=drawColor,line width= 0.4pt,line join=round,line cap=round] (146.93,132.80) circle (  2.50);
\definecolor{drawColor}{RGB}{0,158,115}
\definecolor{fillColor}{RGB}{0,158,115}

\path[draw=drawColor,line width= 0.4pt,line join=round,line cap=round,fill=fillColor] (153.49,148.24) circle (  2.50);

\path[draw=drawColor,line width= 0.4pt,line join=round,line cap=round,fill=fillColor] (160.04,160.23) circle (  2.50);

\path[draw=drawColor,line width= 0.4pt,line join=round,line cap=round,fill=fillColor] (166.59,167.79) circle (  2.50);

\path[draw=drawColor,line width= 0.4pt,line join=round,line cap=round,fill=fillColor] (173.14,169.60) circle (  2.50);

\path[draw=drawColor,line width= 0.4pt,line join=round,line cap=round,fill=fillColor] (179.70,166.07) circle (  2.50);

\path[draw=drawColor,line width= 0.4pt,line join=round,line cap=round,fill=fillColor] (186.25,157.89) circle (  2.50);

\path[draw=drawColor,line width= 0.4pt,line join=round,line cap=round,fill=fillColor] (192.80,146.96) circle (  2.50);
\definecolor{drawColor}{RGB}{213,94,0}

\path[draw=drawColor,line width= 0.4pt,line join=round,line cap=round] (199.35,135.64) circle (  2.50);

\path[draw=drawColor,line width= 0.4pt,line join=round,line cap=round] (205.91,125.79) circle (  2.50);

\path[draw=drawColor,line width= 0.4pt,line join=round,line cap=round] (212.46,118.81) circle (  2.50);

\path[draw=drawColor,line width= 0.4pt,line join=round,line cap=round] (225.56,117.17) circle (  2.50);

\path[draw=drawColor,line width= 0.4pt,line join=round,line cap=round] (238.67,127.91) circle (  2.50);
\definecolor{drawColor}{RGB}{0,158,115}

\path[draw=drawColor,line width= 0.4pt,line join=round,line cap=round,fill=fillColor] (251.77,137.99) circle (  2.50);
\definecolor{drawColor}{RGB}{213,94,0}

\path[draw=drawColor,line width= 0.4pt,line join=round,line cap=round] (264.88,137.99) circle (  2.50);

\path[draw=drawColor,line width= 0.4pt,line join=round,line cap=round] (277.98,137.99) circle (  2.50);

\path[draw=drawColor,line width= 0.4pt,line join=round,line cap=round] (284.54,133.99) circle (  2.50);

\path[draw=drawColor,line width= 0.4pt,line join=round,line cap=round] (291.09,124.18) circle (  2.50);

\path[draw=drawColor,line width= 0.4pt,line join=round,line cap=round] (297.64,112.91) circle (  2.50);

\path[draw=drawColor,line width= 0.4pt,line join=round,line cap=round] (304.19,101.54) circle (  2.50);

\path[draw=drawColor,line width= 0.4pt,line join=round,line cap=round] (310.75, 90.55) circle (  2.50);
\definecolor{fillColor}{RGB}{86,180,233}

\path[fill=fillColor] (314.80, 78.01) --
	(319.80, 78.01) --
	(319.80, 83.00) --
	(314.80, 83.00) --
	cycle;

\path[fill=fillColor] (321.35, 69.49) --
	(326.35, 69.49) --
	(326.35, 74.48) --
	(321.35, 74.48) --
	cycle;

\path[fill=fillColor] (327.91, 62.53) --
	(332.90, 62.53) --
	(332.90, 67.53) --
	(327.91, 67.53) --
	cycle;

\path[fill=fillColor] (334.46, 57.27) --
	(339.45, 57.27) --
	(339.45, 62.27) --
	(334.46, 62.27) --
	cycle;
\definecolor{drawColor}{RGB}{0,0,0}

\path[draw=drawColor,line width= 0.6pt,dash pattern=on 4pt off 4pt ,line join=round] ( 42.21,137.99) -- (356.35,137.99);

\path[draw=drawColor,line width= 0.6pt,dash pattern=on 4pt off 4pt ,line join=round] ( 42.21, 88.51) -- (356.35, 88.51);
\definecolor{drawColor}{RGB}{213,94,0}

\node[text=drawColor,anchor=base,inner sep=0pt, outer sep=0pt, scale=  0.85] at ( 75.53,108.06) {Indecisions};
\definecolor{drawColor}{RGB}{0,158,115}

\node[text=drawColor,anchor=base,inner sep=0pt, outer sep=0pt, scale=  0.85] at ( 75.53,153.04) {Class 1};
\definecolor{drawColor}{RGB}{86,180,233}

\node[text=drawColor,anchor=base,inner sep=0pt, outer sep=0pt, scale=  0.85] at ( 75.53, 63.08) {Class 2};
\end{scope}
\begin{scope}
\path[clip] (  0.00,  0.00) rectangle (361.35,180.67);
\definecolor{drawColor}{RGB}{0,0,0}

\path[draw=drawColor,line width= 0.5pt,line join=round] ( 42.21, 41.95) --
	( 42.21,175.67);
\end{scope}
\begin{scope}
\path[clip] (  0.00,  0.00) rectangle (361.35,180.67);
\definecolor{drawColor}{gray}{0.30}

\node[text=drawColor,anchor=base east,inner sep=0pt, outer sep=0pt, scale=  0.80] at ( 37.71,135.23) {$\tau_\gamma$};

\node[text=drawColor,anchor=base east,inner sep=0pt, outer sep=0pt, scale=  0.80] at ( 37.71, 85.76) {$1-\tau_\gamma$};
\end{scope}
\begin{scope}
\path[clip] (  0.00,  0.00) rectangle (361.35,180.67);
\definecolor{drawColor}{gray}{0.20}

\path[draw=drawColor,line width= 0.5pt,line join=round] ( 39.71,137.99) --
	( 42.21,137.99);

\path[draw=drawColor,line width= 0.5pt,line join=round] ( 39.71, 88.51) --
	( 42.21, 88.51);
\end{scope}
\begin{scope}
\path[clip] (  0.00,  0.00) rectangle (361.35,180.67);
\definecolor{drawColor}{RGB}{0,0,0}

\path[draw=drawColor,line width= 0.5pt,line join=round] ( 42.21, 41.95) --
	(356.35, 41.95);
\end{scope}
\begin{scope}
\path[clip] (  0.00,  0.00) rectangle (361.35,180.67);
\definecolor{drawColor}{RGB}{0,0,0}

\node[text=drawColor,rotate= 90.00,anchor=base,inner sep=0pt, outer sep=0pt, scale=  1.00] at ( 11.89,108.81) {$\eta (x)$};
\end{scope}
\begin{scope}
\path[clip] (  0.00,  0.00) rectangle (361.35,180.67);
\definecolor{fillColor}{RGB}{255,255,255}

\path[fill=fillColor] (117.67,  5.00) rectangle (280.90, 29.45);
\end{scope}
\begin{scope}
\path[clip] (  0.00,  0.00) rectangle (361.35,180.67);
\definecolor{fillColor}{RGB}{255,255,255}

\path[fill=fillColor] (122.67, 10.00) rectangle (137.12, 24.45);
\end{scope}
\begin{scope}
\path[clip] (  0.00,  0.00) rectangle (361.35,180.67);
\definecolor{fillColor}{RGB}{86,180,233}

\path[fill=fillColor] (127.39, 14.73) --
	(132.39, 14.73) --
	(132.39, 19.72) --
	(127.39, 19.72) --
	cycle;
\end{scope}
\begin{scope}
\path[clip] (  0.00,  0.00) rectangle (361.35,180.67);
\definecolor{fillColor}{RGB}{255,255,255}

\path[fill=fillColor] (172.09, 10.00) rectangle (186.55, 24.45);
\end{scope}
\begin{scope}
\path[clip] (  0.00,  0.00) rectangle (361.35,180.67);
\definecolor{drawColor}{RGB}{213,94,0}

\path[draw=drawColor,line width= 0.4pt,line join=round,line cap=round] (179.32, 17.23) circle (  2.50);
\end{scope}
\begin{scope}
\path[clip] (  0.00,  0.00) rectangle (361.35,180.67);
\definecolor{fillColor}{RGB}{255,255,255}

\path[fill=fillColor] (231.47, 10.00) rectangle (245.92, 24.45);
\end{scope}
\begin{scope}
\path[clip] (  0.00,  0.00) rectangle (361.35,180.67);
\definecolor{drawColor}{RGB}{0,158,115}
\definecolor{fillColor}{RGB}{0,158,115}

\path[draw=drawColor,line width= 0.4pt,line join=round,line cap=round,fill=fillColor] (238.70, 17.23) circle (  2.50);
\end{scope}
\begin{scope}
\path[clip] (  0.00,  0.00) rectangle (361.35,180.67);
\definecolor{drawColor}{RGB}{0,0,0}

\node[text=drawColor,anchor=base west,inner sep=0pt, outer sep=0pt, scale=  0.80] at (142.12, 14.47) {Class 2};
\end{scope}
\begin{scope}
\path[clip] (  0.00,  0.00) rectangle (361.35,180.67);
\definecolor{drawColor}{RGB}{0,0,0}

\node[text=drawColor,anchor=base west,inner sep=0pt, outer sep=0pt, scale=  0.80] at (191.55, 14.47) {Indecision};
\end{scope}
\begin{scope}
\path[clip] (  0.00,  0.00) rectangle (361.35,180.67);
\definecolor{drawColor}{RGB}{0,0,0}

\node[text=drawColor,anchor=base west,inner sep=0pt, outer sep=0pt, scale=  0.80] at (250.92, 14.47) {Class 1};
\end{scope}
\end{tikzpicture}}
    \end{minipage}
    \hfill
    \begin{minipage}[t]{0.48\textwidth}
        \centering
        \scalebox{0.65}{%
        \input{Figures/intro_plt_eagle}}
    \end{minipage}
    \captionof{figure}{An example of a binary classification problem that includes indecisions (orange / open circle). In the leftmost figure, the indecisions lie in a region between the two classes: class~1 (green / solid circle) and class~2 (blue / square). A plateau at the threshold $\tau_\gamma$ indicates that some observations may be randomly assigned to either class~1 or the indecision set. In the rightmost figure, the indecisions lie below the threshold $\tau_\gamma$ for the largest class-posterior probability. This demonstrates that the indecision region may not be a simple interval.}
    \label{fig:combined_intro_plots}
\end{figure}
\subsection{Fixed Error Rate}\label{sec:binary_class_fixed_err}

We will now show that understanding the case with a fixed amount of indecisions will allow us to control the accuracy or, similarly, the misclassification rate, of the optimal classifier at any user-specified level. We start with the following result on the properties of the risk function.
\begin{proposition}\label{prop:risk}
    For any $0\leq  \gamma <  1 $, we have 
    $$
    \mathcal{R}(\gamma)= \mathbf{E}_Y\left( Z \; |\; Z < 1 - \tau_{\gamma} \;\; \text{or} \;\;[Z = 1 - \tau_{\gamma}\; \text{and} \; X \notin \mathcal{M}_{\gamma}] \right) 
    ,$$
    where $Z:= \frac{ p_1 f_1 \wedge p_{2} f_{2} }{p_1 f_1 + p_{2} f_{2}}(X) = (\eta\wedge (1-\eta))(X)$.
Moreover, $\gamma \mapsto \mathcal{R}(\gamma)$ is continuous and non-increasing.
\end{proposition}

Because function $\mathcal{R}(\gamma)$ is non-increasing and lower-bounded by $0$, it has a limit as $\gamma \to 1$ that we shall denote $\mathcal{R}^*:=\underset{\gamma \to 1^{-}}{\lim } \mathcal{R}(\gamma) $.
 We note that $\mathcal{R}(\gamma)$ interpolates between the misclassification rate we would get without indecisions and $\mathcal{R}^*$. Thanks to the continuity of $\mathcal{R}$, our result also shows that for any given reachable misclassification level~$\alpha$ above $\mathcal{R}^*$, we can find a~$\gamma^*$ such that $\mathcal{R}(\gamma^*)=\alpha$, and this~$\gamma^*$ is the smallest possible. In other words, for any reachable level of accuracy, we are able to characterize the corresponding minimum number of indecisions. 

\begin{lemma}\label{lem:0risk}
    Suppose that $\mathbf{P}( (\eta \wedge (1-\eta))(X) \leq \varepsilon
    ) >0$  for every $\varepsilon >0$. Then, $\underset{\gamma \to 1^-}{\lim} \tau_{\gamma} = 1$ and $\mathcal{R}^* = 0$.
  \end{lemma}  
    \begin{proof}
        Given any $\varepsilon >0$, there exists a~$\gamma$ such that
    $$
     \mathbf{P}((\eta \wedge (1-\eta))(X) \leq \varepsilon ) > 1 - \gamma \geq  \mathbf{P}((\eta \wedge (1-\eta))(X) < 1 - \tau_{\gamma} ).
    $$
    Consequently, $\tau_{\gamma} \geq 1 - \varepsilon$. Hence,  $
     \mathcal{R}^* \leq  \underset{\gamma \to 1^-}{\lim} 1- \tau_{\gamma} =0
    $ by Proposition~\ref{prop:risk}.  
    \end{proof}
    We conclude that any positive target error level can be reached under the assumption of Lemma~\ref{lem:0risk}. This assumption can be interpreted as follows: there must be regions in which the posterior probability of at least one class is arbitrarily close to one. If definitive predictions from both classes are required, the corresponding condition must hold in both class-dominant regions.

\subsection{Calibration}\label{sec:binary_calib}

\begin{algorithm}[t!]
\caption{Binary Classification with Indecisions}\label{Binary_Class_Calib:alg} 
\noindent\textbf{Input}: Observed $\{(X_i, Y_i): i\in\mathcal{D}\}$, accuracy level  $1-\alpha$. \\
\textbf{\;\;\;\;\; Output}: a selective classification rule $\{\hat Y\in\{0, 1, 2\}\}$ and the corresponding $\hat{\tau}$.

\begin{algorithmic}[1]
    \State Randomly split $\mathcal{D}$ into $\mathcal{D}^{train}$ and $\mathcal{D}^{cal}$.
    \State Train a machine learning model on $\{(X_i, Y_i): i\in \mathcal{D}^{train}\}$. 
    \State Predict the conditional probability $\hat \eta_i$ for all $i\in\mathcal{D}^{cal}$, and set $m=|\mathcal{D}^{cal}|$.
    \State Order $\hat{\tau}_i := 1 - (\hat \eta_i \wedge (1- \hat \eta_i))$ from smallest to largest, $ \hat{\tau}_{(1)}\leq \dots \leq \hat{\tau}_{(m)}$, grouping tied scores together.
    \State For every candidate threshold $\hat{\tau}_{(i)}$, compute the corresponding candidate estimator $\hat{Y}_{(i)}$:
    $$
    \hat{Y}_{(i)}=\mathbf{1}(\hat\eta>\hat\tau_{(i)})+2\,\mathbf{1}(\hat\eta\leq 1 - \hat\tau_{(i)}).
    $$
    \State For every candidate estimator $\hat{Y}_{(i)}$, compute the empirical conditional misclassification error $\hat R_i$ on $\mathcal{D}^{cal}$; discard candidates with no definitive decisions.
    \State Return the first candidate in this ordering for which $\hat R_i\leq\alpha$, thereby minimizing empirical indecision among the feasible candidates.
    \State If no candidate is feasible, report that no empirically feasible rule was found.
\end{algorithmic}

\end{algorithm}

Algorithm~\ref{Binary_Class_Calib:alg} is an empirical calibration heuristic: using the same calibration labels to search over thresholds does not, by itself, provide a finite-sample population-risk guarantee.

Our goal in this section is to present a calibration procedure in the practical setting where the data-generating process is unknown. To do this, we follow the theoretical framework presented in Sections~\ref{sec:binary_class_fixed_indec} and~\ref{sec:binary_class_fixed_err}. This provides a way to empirically calibrate selective classification rules according to our proposed theory. 

We start with the misclassification risk. Given a misclassification error level $\alpha\in[\mathcal{R}^*,\mathcal{R}(0)]$, our goal is to construct a classifier that achieves misclassification level~$\alpha$ using the minimal number of indecisions. From the results in the previous section, we know that there exists an indecision level~$\gamma^*$ such that $\mathcal{R}(\gamma^*) = \alpha$. We aim to construct a classifier~$\hat{Y}$ such that the accuracy of $\hat{Y}$ is at least $1-\alpha$ and the proportion of indecisions is~$\gamma^*$.

Let us define $\gamma_\alpha := \gamma^*$ and $\tau_\alpha:=\tau_{\gamma_\alpha}$ to re-emphasize the fact that the optimal amount of indecisions depends on our desired level of accuracy. We recall that the optimal indecision region is such that $\mathbf{P}\left( \eta(X)\wedge(1-\eta(X))> 1-\tau_\alpha \right) + \mathbf{P}(\mathcal{M}_{\gamma_\alpha})
= \gamma_\alpha.$
Observe that $\tau_\alpha$ corresponds to a quantile of the random variable $\eta(X)\wedge(1-\eta(X)) $ and can be easily computed.
The (conditional) misclassification error of $Y_{\gamma_\alpha}^*$ is ${\mathbf P}_Y( Y_{\gamma_\alpha}^* \neq Y | Y_{\gamma_\alpha}^* \neq 0 ) = \alpha.$ Since we do not have access to $\gamma_\alpha$ explicitly, we need to invert the function $\mathcal{R}(\cdot)$. In order to mimic the optimal classifier $Y_{\gamma_\alpha}^*$, and given an estimator $\hat{\eta}$, we wish to calibrate classifier~$\hat{Y}$ of the form $\hat Y = \mathbf{1}(\hat{\eta}>\hat{\tau}) + 2\,\mathbf{1}(\hat{\eta}\leq 1-\hat{\tau})$, with observations between the two thresholds assigned to the indecision class. We estimate the value of~$\hat{\tau}$ by using a calibration dataset,  as we describe in Algorithm~\ref{Binary_Class_Calib:alg}.

\section{Controlling Type I and Type II Errors: Connection with the Neyman--Pearson Paradigm}\label{sec:gen_hyp_test}
\begin{figure}[t]
    \centering 
    \scalebox{0.8}{
    \input{Figures/np_intro_plt}} \caption{\label{np_intro.plot}\small A comparison of Neyman--Pearson (NP) classification (left) and selective classification with indecisions (right). The NP classifier controls the Type~I error at the target level, at the cost of a larger Type~II error. In contrast, selective classification targets both Type~I and Type~II errors through the introduction of indecisions (yellow shaded region).}
\end{figure}

We now consider the hypothesis-testing problem in which, given a fixed Type~I error probability, the goal is to achieve a desired Type~II error probability, using indecisions if necessary. 

In the absence of indecisions, this framework can be viewed as an alternative to the classical Neyman--Pearson classification paradigm, which prioritizes controlling the Type~I error, often at the expense of the Type~II error \citep{Cannon2002, Scott_05, Rigollet_11, Tong_13}. If the practitioner is willing to incorporate indecision into the decision-making process, it becomes possible to simultaneously target both Type~I and Type~II error rates. We demonstrate the difference between the classical setup and our approach with indecisions in Figure~\ref{np_intro.plot}. 

Our selective classifier $\tilde Y$ can take values in $\{0,1,2\}$, where $0$ denotes an indecision. We define~$A_j$ as the set of all points that our classifier selects for each decision: 
$$
A_j=A_j(\widetilde Y)=\left\{x:\widetilde Y(x)=j\right\},\qquad j=0,1,2,
$$
and denote $\mathbf{P}_i$ as the conditional distribution of $X$ given $Y=i$, for $i \in \{\ 1,2\}$. 

In the context of indecisions, for a given value of $\gamma$, we define the overall probability of making an indecision as $\mathbf{P}(\tilde Y =0)=\mathbf{P}(A_0)$, and the selective Type~I and Type~II error probabilities, respectively, as 
\[
\mathcal P_{\mathrm I}(\widetilde Y)
=\mathbf{P}(\widetilde Y=2\mid Y=1,\widetilde Y\neq0)
=\frac{\mathbf{P}_1(A_2)}{\mathbf{P}_1(A_1)+\mathbf{P}_1(A_2)}
\]
and
\[
\mathcal P_{\mathrm {II}}(\widetilde Y)
=\mathbf{P}(\widetilde Y=1\mid Y=2,\widetilde Y\neq0)
=\frac{\mathbf{P}_2(A_1)}{\mathbf{P}_2(A_1)+\mathbf{P}_2(A_2)}.
\]
When $\gamma = 0$, we recover the usual Type~I and Type~II error probabilities without indecisions. 

Given a Type~I error probability~$\alpha_1$ and an indecision level~$\gamma$, the corresponding optimal Type~II value is given by
$$
\mathcal{P}( \alpha_1,\gamma) = \inf_{\tilde{Y}_{{\gamma}}} \mathcal{P}_{\mathrm{II}}( \tilde{Y}_{{\gamma}}) \nonumber
,$$
where the infimum is taken over all classifiers taking values in $\{0,1,2\}$ such that ${\mathbf P}(\tilde Y_\gamma = 0) = \gamma$ and $\mathcal{P}_{\mathrm{I}}(\tilde{Y}_{{\gamma}}) = \alpha_1 $.
Our goal is to control both the Type~I and Type~II error probabilities at user-specified levels $\alpha_1$ and $\alpha_2$, respectively, using the smallest amount of indecisions.

In some situations, such error control can be achieved without indecisions. However, when the problem is sufficiently difficult, it becomes necessary to identify the hardest to classify observations as indecisions in order to meet our objective.

\subsection{Fixed Level of Indecisions}\label{sec:np_fixed_indec}

We first consider the case where the mass of indecisions is fixed at some level $\gamma \in [0,1)$.
For a fixed level of indecisions $\gamma$, since the regions \(A_0,A_1,A_2\) form a partition and \(\mathbf P(A_0)=\gamma\), this implies that $\mathbf P(A_1)+\mathbf P(A_2)=1-\gamma$.

Our thresholding rule must simultaneously consider both classes and is obtained by placing a mass of indecisions $\gamma$ on the ordered \(\eta\)-scale. Observations below the lower threshold are assigned to class~2, and observations above the higher threshold are assigned to class~1. We then only need to decide where to place this fixed mass of indecisions. We parametrize this location by the marginal mass assigned to the lower decision tail,
\[
r=\mathbf P(A_2)\in[0,1-\gamma].
\]
Equivalently, \(r\) is the mass assigned to the lower decision tail. Once \(r\) is fixed, the upper decision tail must have mass
\[
1-\gamma-r,
\]
because the two decision regions together have total mass \(1-\gamma\). The upper \(\eta\)-tail, where the rule predicts class~1, then has marginal mass
$$
1-\gamma-r.
$$
The corresponding optimal Type~II value is given by
\begin{align}
\mathcal{P}( \alpha_1,\gamma) &= \inf_{\tilde{Y}_{\gamma} } \mathcal{P}_{\mathrm{II}}( \tilde{Y}_{\gamma}) \nonumber \\ 
&= \inf_{\tilde Y_{\gamma}} 
\; {\mathbf P}_{2}(\tilde Y_{\gamma} = 1 \mid \tilde Y_{\gamma} \neq 0 )  \nonumber \\
&= \inf_{\tilde Y_{\gamma}} 
\; \frac{{\mathbf P}_{2}(\tilde Y_{\gamma} = 1)}{\mathbf P_2(\tilde Y\neq0)}, \label{equn:t2e_opt}
\end{align}
where the infimum is taken over all classifiers $\tilde{Y}_{\gamma}$ taking values in $\{0,1,2\}$ such that ${\mathbf P}(\tilde Y = 0) = \gamma$ and the selective Type~I error is $\alpha_1$.

We denote by $ Y^*_{\gamma}$ the optimal classifier achieving the infimum in~\eqref{equn:t2e_opt}. By construction, $Y^*_{\gamma}$ has an overall level of indecisions $\gamma$, a selective Type~I error probability of $\alpha_1$, and the smallest corresponding achievable selective Type~II error probability. 

In the result below, we show that the optimal classifier is of the form
$$
Y^*_{\gamma} \coloneq Y^*_{\gamma}(\alpha_1) =  \mathbf{1}\left(\left\{ \eta > \tau_2  \right\} \cup \mathcal{M}^2_{\alpha_1,\gamma}\right) +  2 \times \mathbf{1}\left(\left\{ \eta \leq \tau_1\right\} \setminus \mathcal{M}^1_{\alpha_1,\gamma} \right),
$$

where the sets
\(\mathcal M^{1}_{\alpha_1,\gamma}\subseteq\{\eta=\tau_1\}\)
and
\(\mathcal M^{2}_{\alpha_1,\gamma}\subseteq\{\eta=\tau_2\}\)
are chosen so that, writing
\[
A_1=\{\eta>\tau_2\}\cup \mathcal M^2_{\alpha_1,\gamma},
\qquad
A_2=\{\eta\le\tau_1\}\setminus \mathcal M^1_{\alpha_1,\gamma},
\]
and
\[
A_0
=
\{\tau_1<\eta<\tau_2\}
\cup
\mathcal M^1_{\alpha_1,\gamma}
\cup
\left(\{\eta=\tau_2\}\setminus \mathcal M^2_{\alpha_1,\gamma}\right),
\]
we have
\[
\mathbf P(A_0)=\gamma
\qquad\text{and}\qquad
\frac{\mathbf P_1(A_2)}
{\mathbf P_1(A_1)+\mathbf P_1(A_2)}
= \alpha_1.
\]
or, equivalently,
\[
(1-\alpha_1)\mathbf P_1(A_2)
=
\alpha_1\mathbf P_1(A_1).
\]

If there are no atoms on the boundary, the Type~I error constraint holds with equality. We then select the optimal thresholds $\tau_1$ and $\tau_2$ by optimizing the placement of the fixed indecision mass $\gamma$ to obtain the smallest selective Type~II error; this corresponds to the feasible indecision set that is as far to the right as possible.

Again, we drop the dependence of \(Y^*_{\gamma}\) on \(\alpha_1\) for simplicity of the presentation.

\begin{theorem}\label{thm:NP:lower_bound}
Classifier \(Y^*_{\gamma}\) is optimal for the selective Type~II risk \(\mathcal{P}(\alpha_1,\gamma)\). Moreover, the optimal risk with overall indecision probability~\(\gamma\) and selective Type~I error control at level~\(\alpha_1\) is given by
\[
\mathcal{P}(\alpha_1,\gamma) = \frac{\mathbf P_{2}(Y^*_{\gamma}=1)} {\mathbf P_{2}(Y^*_{\gamma}\neq0)} = \frac{\int_{\{Y^*_{\gamma}=1\}} f_{2}\,d\mu} {\int_{\{Y^*_{\gamma}\neq0\}} f_{2}\,d\mu}.
\]
\end{theorem}

\subsection{Fixed Error Rate}\label{sec:np_fixed_err_rate}
We now focus on the scenario where a practitioner wishes to control the Type~I and Type~II error probabilities simultaneously at potentially different user-specified levels. We start with a result in the same spirit as Proposition~\ref{prop:risk}, illustrating the connection between the case considered in Section~\ref{sec:np_fixed_indec} and the desired control of these errors. 
\begin{proposition}\label{prop:NP}
    For any $\alpha_1 \in [ 0,1 ]$, the map $\gamma \mapsto \mathcal{P}(\alpha_1, \gamma)$ is continuous and non-increasing.   
\end{proposition}
Since $\mathcal{P}(\alpha_1, \cdot)$ is non-increasing and lower-bounded by~$0$, $\underset{\gamma \to 1^{-}}{\lim } \mathcal{P}(\alpha_1, \gamma)$ exists and shall be denoted $\mathcal{P}^*(\alpha_1)$. Hence, any reachable Type~II error in $(\mathcal{P}^*(\alpha_1),\mathcal{P}(\alpha_1,0)]$ can be obtained using the necessary amount of indecision; the endpoint $\mathcal{P}^*(\alpha_1)$ is included only when that limit is attained for some $\gamma<1$.
\subsection{Calibration}  \label{sec:np_calib}  

\begin{algorithm}[t!]
\caption{Neyman--Pearson Classification with Indecisions}\label{NP_Indec_Calib:alg} 
\hspace*{\algorithmicindent} 

\noindent\textbf{Input}: $\{(X_i,Y_i):i\in\mathcal D^{\mathrm{train}}\}$, $\{(X_i,Y_i):i\in\mathcal D^{\mathrm{cal}}\}$, and Type~I and Type~II targets $\alpha_1,\alpha_2$. \\
\textbf{\;\;\;\;\; Output}: a selective classification rule $\left\{\hat Y\in\{0, 1, 2\}\right\}$ and the corresponding $\hat{\tau}_1$, $\hat{\tau}_2$.

\begin{algorithmic}[1]
    \State Train a machine learning model on $\{(X_i, Y_i): i\in \mathcal{D}^{train}\}$. 
    \State Predict $\hat \eta_i$ for $i\in\mathcal{D}^{cal}$, set $n=|\mathcal D^{cal}|$, and order the scores as $\hat \eta_{(1)}\leq\dots\leq\hat \eta_{(n)}$, grouping ties together. Set $\Gamma=\{k/n:k=0,\dots,n-1\}$.
    \State For each $\gamma_k\in\Gamma$, consider the score-cutpoint pairs that form a candidate indecision block of mass $\gamma_k$ (with tied scores kept together). Discard candidates having a zero class-specific decision count.
    \State Among candidates whose empirical selective Type~I error is at most $\alpha_1$, choose the rightmost feasible block and record its thresholds as $\hat\tau_1(k)$ and $\hat\tau_2(k)$. If none is feasible, discard $\gamma_k$.
    \State Compute the corresponding candidate estimator $\hat{Y}_{\gamma_k}$:
    $$
    \hat{Y}_{\gamma_k}=\mathbf{1}(\hat\eta>\hat\tau_2(k))+2\,\mathbf{1}(\hat\eta\leq\hat\tau_1(k)).
    $$
    \State Estimate the empirical selective Type~II error of $\hat{Y}_{\gamma_k}$.
    \State Among the retained candidates with empirical Type~II error at most $\alpha_2$, choose the smallest $\gamma_k$ and denote it by $\hat\gamma$.
    \State If no candidate remains, report that no empirically feasible rule was found; otherwise return $\hat{Y}_{\hat\gamma}$ and $\hat\tau_1(\hat\gamma),\hat\tau_2(\hat\gamma)$.
\end{algorithmic}
\end{algorithm}

Like Algorithm~\ref{Binary_Class_Calib:alg}, Algorithm~\ref{NP_Indec_Calib:alg} provides empirical calibration only; a population guarantee requires a separate finite-sample certification argument.

We are now equipped to use the results in Sections~\ref{sec:np_fixed_indec} and~\ref{sec:np_fixed_err_rate} to calibrate thresholds for Type~I and Type~II error targets. In contrast to the binary classification setting in Section~\ref{sec:general}, calibration is more challenging because it requires class-conditional ratio estimation, whose denominators may be small. The optimally attainable Type~II value may decrease as the abstention budget grows, but arbitrary added indecision need not lower the Type~II error of a fixed rule.

For given levels $\alpha_1 \in [0,1]$ and $\alpha_2 \in [\mathcal{P}^*(\alpha_1),\mathcal{P}(\alpha_1,0)]$ of Type~I and Type~II error, respectively, there exists a $\gamma^*$ such that $\mathcal{P}(\alpha_1, \gamma^*)=\alpha_2$ when the target is attainable. Our goal is to find a classifier $\hat{Y}$ whose empirical Type~I and Type~II errors are at most~$\alpha_1$ and~$\alpha_2$, respectively, while approximating the minimum indecision mass $\gamma^*$. Our suggested estimator~$\hat Y$ is of the form $\hat{Y} =  \mathbf{1}\left(\left\{ \hat \eta > \hat{\tau}_2  \right\}\right) +  2 \times \mathbf{1}\left(\left\{ \hat \eta \leq \hat{\tau}_1\right\} \right).$

Ideally, we would like the optimal indecision region to be given by $\mathbf{P}\left( \hat{\tau}_1 < \hat{\eta} < \hat{\tau}_2\right) = \gamma^*$.
Because we do not have access to $\gamma^*$ explicitly, we need to invert the function $ \mathcal{P}(\alpha_1, \cdot)$.
We estimate the values $\hat{\tau}_1,\hat{\tau}_2$ using a calibration dataset, as we describe in Algorithm \ref{NP_Indec_Calib:alg}.

\section{Theory}\label{sec:theory}

In this section, we further investigate the theoretical properties of our selective classification framework. The first subsection studies an exponent-level phase transition in the risk of our classifier and sheds light on an all-or-nothing phenomenon. The second considers plug-in rules in which the true conditional probability function $\eta$ is replaced with a learned function $\hat\eta$.

\subsection{Explicit indecisions for the Gaussian Mixture Model: A sharp phase transition}\label{sec:gaussian}

Selective classification is especially attractive when a small indecision mass yields a substantial improvement in accuracy, as illustrated in Figure~\ref{accurary_gamma.plot}. Understanding this regime is important for practical adoption: difficult tasks may defer too many cases to human reviewers, whereas easy tasks may require no indecisions because standard classifiers already meet the target. The greatest benefit therefore arises in an intermediate regime, where standard models fall short of the target but a modest indecision rate is sufficient. In this section, we analyze the exponent-level transition in the selective risk when the Bayes classifier cannot meet the desired accuracy.

This section is devoted to the asymptotic behavior of the optimal amount of indecisions as the risk gets smaller. We focus on the symmetric Gaussian mixture model. In particular, we assume that $p_1 = p_{2} =1/2$ and that $f_1(x) = \frac{1}{\sqrt{2\pi}} \exp(-(x-\Delta)^2/2) = f_{2}(-x)$ for some separation $\Delta>0$. In this case, the likelihood ratio is monotone. We describe the benefits of the monotone likelihood-ratio (MLR) property further in Section~\ref{sec:mlr}. For a given misclassification rate $\delta \to 0$, we are interested in the asymptotic behavior of $\gamma_\delta$ as a function of the separation $\Delta$. Naturally, we expect $\gamma_\delta$ to be non-increasing in~$\Delta$. We assume that $\delta\to
0$ and let parameter~$\Delta$ depend on~$\delta$, omitting subscript~$\delta$ whenever no ambiguity arises. 

The asymptotic property we study here is \textit{$\delta$-consistency}, which is inspired by the consistency in classification of \cite{minsker2021minimax} or, similarly, exact recovery in Gaussian mixtures as defined in \cite{ndaoud2022sharp}. We establish an exponent-level characterization of the phase transition for $\delta$-consistency.
\begin{definition}\label{def:consistency} Let $(\gamma_\delta)$ be indexed by a sequence of positive values $\delta\downarrow0$.

\begin{itemize}
\item We say that \emph{$\delta$-consistency is impossible} for
$(\gamma_\delta)$ if
$$
\label{af2} \liminf_{\delta \to 0}\mathcal{R}(\gamma_\delta)/\delta >1.
$$

\item We say that \emph{$\delta$-consistency is possible} for $(\gamma_\delta)$
 if there exists a classifier  $\hat Y_{\gamma_\delta}:=\hat Y(\gamma_\delta, \cdot) $, such that $\mathbf{P}(\hat Y_{\gamma_\delta} = 0) = \gamma_\delta$ for all $\delta$, and 
$$
\label{af1}  \limsup_{\delta\to0}  {\mathbf P}_Y(\hat Y_{\gamma_\delta}  \neq Y | \hat Y_{\gamma_\delta}  \neq 0 ) /\delta \leq 1.
$$
In this case, we say that $\hat Y$ achieves $\delta$-consistency.

\end{itemize}
\end{definition}

In order to derive the phase transition of interest, let us first recall the equations that relate~$\gamma_\delta$ to~$\Delta_\delta$ and~$\delta$. For a misclassification level~$\delta$, we have, under the MLR property (Section \ref{sec:mlr}), that $\mathbf{P}(\xi \geq \Delta_\delta + t_\delta) = (1-\gamma_\delta) \delta,$ where $\xi$ is a standard normal random variable, and $t_\delta$ is a threshold that can be related to $\tau_\delta$. Moreover, the indecision level is given by $\mathbf{P}(\xi \geq \Delta_\delta - t_\delta) - \mathbf{P}(\xi \geq \Delta_\delta + t_\delta) = \gamma_\delta.$
Since there is a one-to-one correspondence between $\delta$ and $\gamma_\delta$, it is easy to see that the same holds for $t_\delta$ as well. Our proof strategy works as follows. For a given $t\geq 0$:
$$
\frac{\mathbf{P}(\xi \geq \Delta + t)}{\mathbf{P}(\xi \geq \Delta + t) + \mathbf{P}(\xi \geq t - \Delta )} \leq 
\delta \quad \text{ if and only if } \quad \mathbf{P}(\xi \geq \Delta - t) - \mathbf{P}(\xi \geq \Delta + t) \geq\gamma_\delta.
$$

We use the following parameterizations for $\Delta_\delta$ and~$\gamma_\delta$: $\Delta_\delta = c\sqrt{2\log(1/\delta)}$ for $0<c<1$ and $m>0$, and
$$
\gamma_\delta = \left\{
    \begin{array}{ll}
        1 - \delta^{m}& \mbox{ if } 0<c < 1/2, \\
        \delta^{m} & \mbox{ if }  1/2<c < 1.
    \end{array}
\right.
$$
We also define $m^*(c)$ such that
\begin{equation}\label{eq:m(c)}
   m^*(c) = \left\{
    \begin{array}{ll}
        (c-1/(4c))^2& \mbox{ if } 0<c < 1/2, \\
        (2c-1)^2 & \mbox{ if }  1/2<c < 1.
    \end{array}
\right. 
\end{equation}
The next result
describes a ``phase transition'' for $\gamma_{\delta}$ for the problem of $\delta$-consistency.

\begin{theorem}\label{thm:asymp}
For any $\varepsilon>0$ and $1/2<c<1$:
\begin{itemize}
\item[(i)] Let $m \leq  m^*(c). $
Then, the classifier $Y^*_{\gamma_\delta}$ defined in \eqref{eq:oracle} achieves $\delta$-consistency.
\item[(ii)] Moreover, if $m \geq (1+\varepsilon) m^*(c)$, then $\delta$-consistency is impossible.
\end{itemize}
For any $\varepsilon>0$ and $0<c<1/2$:
\begin{itemize}
\item[(i)] Let $m \geq  m^*(c).$
Then, the classifier $Y^*_{\gamma_\delta}$ defined in \eqref{eq:oracle} achieves $\delta$-consistency.
\item[(ii)] Moreover, if $m \leq (1-\varepsilon) m^*(c)$, then $\delta$-consistency is impossible.
\end{itemize}
\end{theorem}
\noindent Theorem~\ref{thm:asymp} gives the following sufficient regimes for $\delta$-consistency; the corresponding impossibility regimes hold with the $\varepsilon$-gaps stated in the theorem:
\begin{align}
\label{phase}
m \leq m^*(c), &\quad \text{for } 1/2 < c < 1, \\
\label{phase2}
m \geq m^*(c), &\quad \text{for } 0 < c < 1/2.
\end{align}

It is worth noting that, while in the classical setup (without indecisions) we need $c\geq1$ to achieve $\delta$-consistency, we require almost no indecisions when $c>1/2$, as $\delta^{(2c-1)^2}=o(1)$. We also observe an all-or-nothing phenomenon: $\gamma_\delta$ tends to $0$ or $1$ depending on whether $c$ is greater or smaller than $1/2$.

\subsection{Plug-in rules}
\label{sec.plugin}

In this section, we provide an oracle plug-in analysis motivated by Algorithm~\ref{Binary_Class_Calib:alg}.
For simplicity, we assume an infinite calibration set, allowing us to focus on the error introduced by estimating~ $\eta$ from a finite training sample.
Analogous results for hypothesis testing (Algorithm \ref{NP_Indec_Calib:alg}) can be derived in the same manner and are therefore omitted to avoid redundancy.

In what follows, we replace classification probability function~$\eta$ with a learned  function~$\widehat{\eta}$. Given an indecision level, we quantify the loss in the accuracy due to the estimation of~$\eta$. In addition, we also investigate what happens to the indecision level if we calibrate the $\widehat{\eta}$-based method to achieve a pre-specified level of accuracy. 

Similar results have been established in the literature. In particular, \cite{Denis_20} derive results on the asymptotic performance of the plug-in classifier in the general setting of our Theorem~\ref{lem.fixed.indec}. However, they assume that $\eta(X)$ has a continuous distribution near the decision threshold, while we allow for a point mass on the decision boundary. Furthermore, \cite{Lei_14} establishes results that are similar to our Theorem~\ref{lem.fixed.acc}. However, while we fix the conditional misclassification error given that a decision has been made, \cite{Lei_14} focuses on the unconditional accuracy. Moreover, like \cite{Denis_20}, he also assumes that $\eta(X)$ has a continuous distribution near the decision threshold.

The accuracy and the level of indecisions are linked through the choice of the threshold~$\tau$. Our analysis is more challenging than the one in \cite{Herbei_Wegkamp_06}, because it relies on ensuring that $\hat{\tau}$ and $\tau$ are relatively close, which requires stronger assumptions than the usual margin condition.  
For the rest of this section, and for the sake of simplicity, we will assume that all levels of misclassification $\alpha>0$ can be reached using our framework. In what follows we shall only provide theoretical results for the risk of the plug-in  classifier in the classification setting. Similar results can be derived for testing as well. 

\subsubsection{Fixing the probability of an indecision} 

Let $\widehat{Y}_{\gamma}$ be the plug-in classifier for the indecision level~$\gamma$, i.e., 
$$
\widehat{Y}_{\gamma}(X)=1\times\mathbf{1}\{ \widehat{\eta}(X)>\widehat{\tau}_{\gamma}\}+2\times\mathbf{1}\{ \widehat{\eta}(X)<1-\widehat{\tau}_{\gamma}\},
$$
where $\widehat{\tau}_{\gamma}$ is chosen so that $\mathbf{P}\Big( \widehat{\tau}_{\gamma}\ge\widehat{\eta}(X)\ge 1-\widehat{\tau}_{\gamma}\,|\,\widehat{\eta}\Big)=\gamma.$ Note that we may need to use only a subset of $\{\widehat{\eta}(X)=\widehat{\tau}_{\gamma}\}\cup \{\widehat{\eta}(X)=1-\widehat{\tau}_{\gamma}\}$ rather than the full set to get the exact equality above, same as we did for $Y^*_{\gamma}$.

We will write $\tau^*_{\gamma}$ for the corresponding threshold for $Y^*_{\gamma}$. 
Given a classifier $\tilde{Y}$, we let $R(\tilde{Y})=\mathbf{P}( \tilde{Y}\ne Y | \tilde{Y}\ne 0)$. 
To simplify expressions, we write~$\eta$ and~$\widehat{\eta}$ instead of~$\eta(X)$ and~$\widehat{\eta}(X)$, respectively. 

\smallskip

\begin{lemma}
\label{lem.risk.expr}
For all $\gamma\in[0,1)$,

\vspace{10pt}

$R(\widehat{Y}_{\gamma})-R(Y^*_{\gamma}) \, = \,\tfrac{1}{1-\gamma}\mathbf{E}|\tau^*_{\gamma}-\eta|\big(\mathbf{1}\{Y^*_{\gamma}=1,\widehat{Y}_{\gamma}\ne Y^*_{\gamma}\}+\mathbf{1}\{\widehat{Y}_{\gamma}=1,\widehat{Y}_{\gamma}\ne Y^*_{\gamma}\}\big)\,+$\\

$~~~~~~~~~~~~~~~~~~~~~~~~~~~~~~~\tfrac{1}{1-\gamma}\mathbf{E}|1-\tau^*_{\gamma}-\eta|\big(\mathbf{1}\{Y^*_{\gamma}=2,\widehat{Y}_{\gamma}\ne Y^*_{\gamma}\}+\mathbf{1}\{\widehat{Y}_{\gamma}=2,\widehat{Y}_{\gamma}\ne Y^*_{\gamma}\}\big)$.

\end{lemma}

\begin{remark}\label{rmk:upper}
    We can bound the above expression as follows:
    $$
R(\widehat{Y}_{\gamma})-R(Y^*_{\gamma}) \le 2\mathbf{E}\big((|\tau^*_{\gamma}-\eta| \vee  |1-\tau^*_{\gamma}-\eta|)\mathbf{1}\{\widehat{Y}_{\gamma}\ne Y^*_{\gamma}\} \big|
\widehat{Y}_{\gamma} \neq 0 \;\textup{or}\; Y^*_{\gamma} \neq 0\big).
$$
\end{remark}

Remark \ref{rmk:upper} 
implies that if~$\eta$ does not
have too much mass around the optimal thresholds $\tau^*_{\gamma}$ and $1-\tau^*_{\gamma}$, and~$\hat{\eta}$ is close to $\eta$, 
then we can expect consistency of the plug-in approach. 

For the next result, we let $\eta_{\text{max}}:=\eta\vee(1-\eta)$ and focus on the standard setting where we can bound the probability that~$\eta_{\text{max}}$ lies within~$\phi$ of $\tau^*_{\gamma}$ by some nonnegative power of $\phi$. More specifically, we assume
\begin{equation}
    \mathbf{P}\big(| \eta_{\text{max}} - \tau^*_{\gamma} | \leq \phi \big)\lesssim \phi^{\beta} \quad\text{and}\quad 
    \mathbf{P}\big(\tau^*_{\gamma}< \eta_{\text{max}} \le \tau^*_{\gamma} +\phi \big) \wedge  \mathbf{P}\big(\tau^*_{\gamma}-\phi\le \eta_{\text{max}} < \tau^*_{\gamma}  \big) \gtrsim\phi^{\beta'}\label{lbnds_eta_prob01},
\end{equation}
for some $\beta' \geq \beta\ge0$ and all sufficiently small positive~$\phi$.

\begin{theorem}
\label{lem.fixed.indec}
Suppose that~(\ref{lbnds_eta_prob01}) holds for $0<\phi\le3\phi^*_\gamma$ with $\beta' \geq \beta\ge0$ and $\phi^*_{\gamma}>0$. Then,
\begin{equation}
\label{cl1.lem.fixed.indec}
    R(\widehat{Y}_{\gamma})-R(Y^*_{\gamma}) \lesssim \inf_{0<\phi\le\phi^*_\gamma} \left\{\tfrac{1}{1-\gamma}\Big[\mathbf{P}\big(|\widehat{\eta}-\eta|>\phi\big)+\phi^{1+\beta}\Big]+\Big[\phi^{1-\beta'}\mathbf{P}\big(|\widehat{\eta}-\eta|>\phi\big) \wedge \phi\Big]\right\}.
\end{equation}
In particular, if $\beta' \le1$, then
$$
R(\widehat{Y}_{\gamma})-R(Y^*_{\gamma}) \lesssim \tfrac{1}{1-\gamma}\inf_{0<\phi\le\phi^*_\gamma}\left\{\mathbf{P}\big(|\widehat{\eta}-\eta|>\phi\big)+\phi^{1+\beta}\right\}.
$$
\end{theorem}

\begin{remark}
    In the statement of Theorem~\ref{lem.fixed.indec}, we can replace $\tfrac{1}{1-\gamma}\mathbf{P}\big(|\widehat{\eta}-\eta|>\phi\big)$ with $\mathbf{P}\big(|\widehat{\eta}-\eta|>\phi \,|\, \eta_{\text{max}}>\tau^*_{\gamma}\big)+\mathbf{P}\big(|\widehat{\eta}-\eta|>\phi \,|\, \widehat{\eta}_{\text{max}}>\widehat{\tau}_{\gamma}\big)$, where $\widehat{\eta}_{\text{max}}:=\widehat{\eta}\vee(1-\widehat{\eta})$. That is, we only need $\hat{\eta}$ to be close to $\eta$ within the region of decisions. 
    For example, this can be easily achieved if we have good  control over the uniform bound $\|\widehat{\eta}-\eta\|_\infty$. We can also replace the term $\frac{\phi^{1+\beta}}{1-\gamma}$ by $\phi^{1+\beta}$ if we assume that the margin condition \eqref{lbnds_eta_prob01} holds conditionally on being in the region of decisions.
\end{remark}

Note that we will have a good estimator $\hat{\eta}$ of $\eta$ as long as $\eta$ is sufficiently smooth. When $\beta' \leq 1$, our result is similar to the corresponding one in \cite{Herbei_Wegkamp_06}, which covers the setting without indecisions. Unlike that setting, we face the additional challenge of controlling the distance between thresholds~$\hat{\tau}$ and~$\tau$. The lower bound in condition \eqref{lbnds_eta_prob01} provides that control. Consequently, when $\phi \approx 1/\sqrt{n}$, where~$n$ is the training sample size, we can recover fast rates when $\beta = \beta' =1$, which is typically the case for atom-less distributions. Returning to~\eqref{cl1.lem.fixed.indec} and taking $\phi \approx 1/\sqrt{n}$,  we recover the slow rates without margin assumptions.

\smallskip

\subsubsection{Fixing the misclassification level}
We will use $R_{\widehat{\eta}}(\widehat{Y})$ to denote the conditional risk of the classifier~$\widehat{Y}_{\gamma}$ given $\widehat{\eta}$.
Let~$\gamma$ be a fixed indecision level and let $R^*:=R(Y^*_{\gamma})$.
Here, we analyze the plug-in classifier corresponding to the misclassification level~$R^*$. The classifier we consider is of the form
$1\times\mathbf{1}\{ \widehat{\eta}>\widehat{\tau}\}+2\times\mathbf{1}\{ \widehat{\eta}<1-\widehat{\tau}\}$. Its indecision level and threshold are respectively given by $\widehat{\gamma}:=\inf\{u:  R_{\widehat{\eta}}(\widehat{Y}_{u})\le R^* \}$ and  $\widehat{\tau}:=\inf\{t: \mathbf{P}( {t}\ge\widehat{\eta}\ge 1-{t}\,|\,\widehat{\eta})\ge\widehat{\gamma}\}$. 

\begin{theorem}
\label{lem.fixed.acc}
Suppose that $\gamma<1$ is a fixed indecision level and condition (\ref{lbnds_eta_prob01}) holds for $0<\phi\le2\phi^*_{\gamma}$ with $\beta' \geq \beta\ge0$ and $\phi^*_{\gamma}>0$. Then, there exists a positive universal constant~$c_1$ such that
$$
\mathbf{P}\big(\widehat{\tau}-\tau^*_{\gamma}>\phi\big)\lesssim
\frac{\mathbf{P}\big(|\widehat{\eta}-\eta|>c_1\phi^{1+2\beta'-2\beta}\big)}{\phi^{1+2\beta'-\beta}},
$$
for $0<\phi\le\phi^*_{\gamma}$. Moreover, we also have that
$$
\mathbf{E}(\hat{\gamma}) - \gamma \lesssim
 \inf_{0<\phi\le\phi^*_{\gamma}}\left\{ \phi^\beta + \frac{\mathbf{P}\big(|\widehat{\eta}-\eta|>c_1\phi^{1+2\beta'-2\beta}\big)}{\phi^{1+2\beta'-\beta}}\right\}.
$$
\end{theorem}

In the case $\beta=\beta'=1$, we can expect to recover slow rates of classification, while in general consistency is not guaranteed, especially if~$\eta_{\text{max}}$ has some mass around~$\tau^*_\gamma$.

\section{Extensions}\label{appendix:sec:ext}

First, we demonstrate that selective classifiers described in Sections \ref{sec:general} and \ref{sec:gen_hyp_test} can completely avoid estimation of the regression function, through the monotone likelihood ratio (MLR) property. Second, we consider the multiclass classification case and derive the minimax rules for this setting. 

\subsection{Adaptation under the MLR property}\label{sec:mlr}

In both the binary and Neyman--Pearson classification settings, we can avoid estimating the conditional probability function $\eta(\cdot)$ in \eqref{cond_dens.equn}.

For simplicity, we assume for the remainder of this section that the distribution of~$X$ does not have any atoms, i.e., $\mathbf{P}(X=t)=0 \;\forall t\in \mathbb{R}$. 
While the approach discussed above allows us to calibrate the optimal procedure with indecisions, it relies heavily on prior knowledge of the likelihood ratio $p_1f_1/(p_2f_2)$. 
In this section, we demonstrate how to achieve adaptation under the monotone likelihood-ratio (MLR) property, which is defined as follows:

{\it The random variable~$X$ takes values in a subset of~$\mathbb{R}$, the densities~$f_1$ and~$f_{2}$ have the same support, and $\frac{f_{2}}{f_{1}}(\cdot)$ is an increasing function on the support of the densities.}

This property covers a large class of exponential-family distributions. For example, it is satisfied for location models with a log-concave density. It is also satisfied for the chi-square location model, where $f_{2}$ and $f_{1}$ are, respectively, standard chi-square and noncentral chi-square densities. We refer the reader to \cite{butucea2023variable} for more details about the MLR property.

More precisely, under the MLR property we can calibrate the oracle procedure based on observations of $X$ and without prior knowledge of $f_{1}$ or $f_2$. For the Neyman--Pearson testing problem, the optimal procedure under MLR is given by
$$
Y_{\gamma}^* = \mathbf{1}(X \leq \tau_2 ) + 2 \times \mathbf{1}(X \geq \tau_1 ),
$$
where $\tau_1,\tau_2$ are such that 
$$
\frac{\mathbf{P}_1(X\geq \tau_1)}
{\mathbf{P}_1(X\leq\tau_2)+\mathbf{P}_1(X\geq\tau_1)}=\alpha_1
\quad\text{and}\quad
\mathbf{P}\left( \tau_2\leq X \leq \tau_1 \right) = \gamma.
$$
The Type~II error of $Y_{\gamma}^*$ is given by
$$
 \mathcal{P}(\alpha_1,\gamma) =
 \frac{{\mathbf P}_{2}(X \leq \tau_2)}
 {{\mathbf P}_{2}(X\leq\tau_2)+{\mathbf P}_{2}(X\geq\tau_1)}.
$$
Recall that our goal is to calibrate $\gamma$ so that $\mathcal{P}(\alpha_1,\gamma) = \alpha_2$.
\begin{remark}
    Observe that $\gamma > 0$ (i.e., $\tau_2 < \tau_1$) is needed when $F^{-1}_2(\alpha_2) < F^{-1}_1(1-\alpha_1)$. In other words, indecisions are needed if the power of the nonselective NP test is below the target $1 - \alpha_2$.
\end{remark}

Given a calibration set of i.i.d. $X_i$ and the corresponding labels, we can compute the above quantiles empirically and repeat the steps described above. It is interesting to note that under the MLR property the indecision set is an interval. 

The case with accuracy is slightly more challenging, as the constraint 
$$
\mathbf{P}\left( 1 - \tau_{\gamma }  \leq  \eta(X)\leq \tau_{\gamma } \right) = {\gamma }
$$
does not necessarily translate into a symmetric interval for $X$. This can be dealt with if we further assume that $\log\left(\frac{p_1f_1}{p_{2}f_{2}} \right)(\cdot)$ is an odd function. In particular, this is the case under mixtures of symmetric distributions such that $p_{2}f_{2}(x) = p_1f_1(-x)$. In that case, the optimal procedure becomes
$$
Y_{\gamma }^* = 2 \times \mathbf{1}(X\geq \tau_{\gamma }  ) + \mathbf{1}(X \leq - \tau_{\gamma }  ),
$$
where $\tau_{\gamma }\in [0,\infty)$ is such that 
$$
1 - 2\mathbf{P}(X \geq \tau_{\gamma })=\mathbf{P}\left( -\tau_{\gamma }  \leq  X \leq \tau_{\gamma } \right) = {\gamma }.
$$
Our goal here is to calibrate $\gamma$ so that the misclassification rate of $Y_{\gamma }^*$ satisfies
$$
{\mathbf P}_Y( Y_{\gamma }^* \neq Y | Y_{\gamma }^* \neq 0 ) := \frac{p_1 \mathbf{P}_1\left( X \geq  \tau_{\gamma } \right)+ p_{2}\mathbf{P}_{2}\left( X\leq -\tau_{\gamma } \right)}{1 - {\gamma }} = \alpha.
$$
Again, using a calibration set, we can estimate the above quantiles and recover the optimal classifier under indecisions.

\subsection{Multiclass Classification}\label{sec:multi}
We now focus on the multiclass case. 
Assume that we observe a random variable $X$ on a measurable space $(\mathcal{X},\mathcal{U})$ such that $X$ is distributed according to a mixture model, where with probability $p_i$ its probability measure is given by $P_{i}$ for $i=1,\dots,K$, and~$K$ is the number of classes. We assume that $P_{i}\neq P_{j}$ for any $i\neq j$. Let~$f_{i}$ be the density of $P_{i}$ with respect to some dominating measure that we will further denote by $\mu$. Denote by $Y$ the labeling quantity such that $Y=i$ if the distribution of $X$ is $P_{i}$. We are interested in the problem of recovering the label $Y$. 

As estimators of $Y$, we consider all measurable functions $\widehat Y=\widehat Y(X)$ of $X$  taking values in $\{0,1,\dots,K\}$, where we allow for indecisions. 
Such estimators will be called {\it classifiers}. 
The performance of a classifier $\widehat Y$ is measured by its expected risk ${\mathbf P}_Y (\widehat Y \neq Y\mid \widehat Y \neq 0 )$. We denote by ${\mathbf E}_Y$ the expectation with respect to probability measure ${\mathbf P}_Y$ of $X$ with label $Y$. Observe that
$$
{\mathbf P}_Y (\widehat Y \neq Y \mid \widehat Y \neq 0 ) = \frac{\sum_{i=1}^K p_i {\mathbf P}_i (\widehat Y \notin \{i,0\})}{\mathbf{P}(\widehat Y\neq 0)} = 1 - \frac{\sum_{i=1}^K p_i {\mathbf P}_i (\widehat Y = i)}{\mathbf{P}(\widehat Y\neq 0)}.
$$
For an indecision level $\gamma$, the optimal selective risk $\mathcal{R}$ is given by
$$
\mathcal{R}(\gamma) := \inf_{\tilde Y}  {\mathbf P}_Y(\tilde Y \neq Y | \tilde Y \neq 0 ),
$$
where $\inf_{\tilde Y}$ denotes the infimum over all classifiers taking values in $\{0,1,\dots,K\}$ such that ${\mathbf P}(\tilde Y = 0) = \gamma$. Define the oracle classifier by
\begin{equation}\label{eq:multi}
    Y_\gamma^* = \mathbf{1}\left(\underset{i}{\max}(p_if_i) \geq \tau_\gamma \sum_i p_i f_i\right)\underset{i}{\arg\max(p_if_i)},
\end{equation}
where $\tau_\gamma \in [1/K,1]$ is such that 
$$
\mathbf{P}\left( \underset{i}{\max}(p_if_i) \leq \tau_\gamma \sum_i p_i f_i \right) = \gamma.
$$

\begin{theorem}\label{thm:multi:lower_bound}
   The classifier $Y_\gamma^*$, defined in \eqref{eq:multi}, is optimal for the risk $\mathcal{R}(\gamma)$. Moreover, we have
   $$ 
   \mathcal{R}(\gamma) = {\mathbf P}_Y( Y_\gamma^* \neq Y |  Y_\gamma^* \neq 0 ) = 1 - \frac{\int_{\{Y_\gamma^* \neq 0\}}\underset{i}{\max}( p_if_i)d\mu }{1-\gamma}.
   $$
\end{theorem}

\begin{remark}
    We can use plug-in scores to calibrate the procedure as we did in Sections~\ref{sec:binary_calib}, \ref{sec:np_calib}, and \ref{sec.plugin}. Note that, for a given indecision level $\gamma$, the procedure does not require labels for calibration and can therefore be calibrated in an unsupervised fashion.
\end{remark}

\section{Simulations}\label{sec:simu}
We illustrate the theory established in this paper through three simulation studies. The first demonstrates that we can empirically recover the phase transition in Section~\ref{sec:gaussian} under the two-component Gaussian mixture model. The second demonstrates the hypothesis-testing procedure detailed in Section~\ref{sec:np_fixed_err_rate}. In these simulations, adding an indecision region often lowers the Type~II error when the Type~I target is met; this is an empirical observation about the fitted procedures, not a general monotonicity guarantee.

\subsection{Phase Transition: Gaussian Mixture Model}\label{sim:phase_transition}

\begin{figure}[!t]  
    \centering
    \includegraphics[scale=0.8]{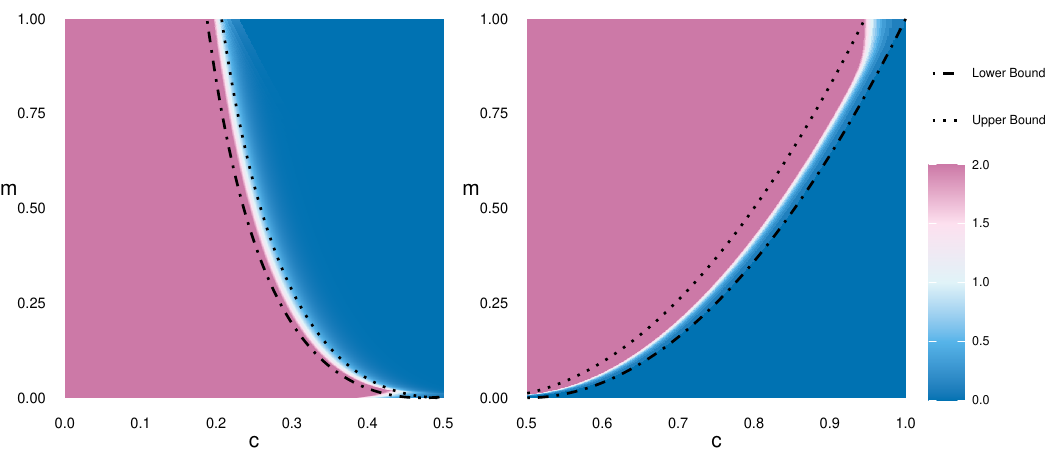} 
    \caption{Computation of $\mathcal{R}(\gamma)/\delta$ with $\delta=10^{-7}$ (left) and $\delta=10^{-15}$ (right). The lower bound corresponds to the curve $m_*(c)$, while the upper bound corresponds to $m^*(c)$.}
    \label{fig:transition}
\end{figure}

In the symmetric Gaussian mixture model, given a vanishing misclassification level~$\delta$, we compare the theoretical and empirical values of~$\gamma_{\delta}$.  
We recall the parameters
 $c$ and $m$ such that 
$\Delta(c) = c\sqrt{2\log(1/\delta)}$
and $\gamma_{\delta}(m) =(1 - \delta^{m})\mathbf{1}\{0<c < 1/2\} + \delta^{m}\mathbf{1}\{1/2<c < 1\}$, where we have made the dependence of $\gamma_{\delta}$ on~$m$ explicit.

According to Theorem \ref{thm:asymp}, and using the above parameterization, $\delta$-consistency is possible whenever
$$
m \geq m^*(c) \quad \text{if} \quad 0<c < 1/2 \qquad\text{and}\qquad  
       m \leq m^*(c) \quad\text{if}\quad  1/2<c < 1.
$$
Based on the proof of Theorem \ref{thm:asymp}, it also holds that $\delta$-consistency is impossible whenever
$$
m \leq m_*(c) \quad \text{if} \quad 0<c < 1/2 \qquad\text{and}\qquad  
       m \geq m_*(c) \quad\text{if}\quad  1/2<c < 1,
$$
where the lower bound $m_*(c)$ is given by
\begin{equation}
   m_*(c) = \left\{
    \begin{array}{ll}
        (c-(1-\varepsilon)/(4c))^2& \mbox{ if } 0<c < 1/2, \\
        (2c-1+\varepsilon)^2 & \mbox{ if }  1/2<c < 1,
    \end{array}
\right.
\end{equation}
for $ \varepsilon= 1/2\log(4 \pi  \log(1/\gamma_{\delta}))/\log(1/\gamma_{\delta})$. Observe that $m_*(c)\to m^*(c)$ as $\delta \to 0$.

In what follows, we
fix $\delta=10^{-15}$ for $c>1/2$, and $\delta=10^{-7}$ for $c<1/2$. 
Our simulation setup is defined as follows. We set $c$ on a uniform grid
of $1{,}000$ points delimited by $0$ and $1$. Similarly, we set~$m$ on a uniform grid
of $1{,}000$ points delimited by~$0$ and~$1$. For each combination of values of $c$
and~$m$, we first find $t_\gamma$ using a grid search such that $\mathbf{P}(\xi \geq \Delta(c) - t_\gamma) -  \mathbf{P}(\xi \geq \Delta(c) + t_\gamma) =  \gamma_{\delta}$, where~$\xi$ is a standard normal random variable. 
Next, we compute $\mathcal{R}(\gamma_{\delta})= \mathbf{P}(\xi \geq \Delta(c) + t_\gamma) / (1-\gamma_{\delta}).$
To improve interpretability, values of $\mathcal{R}(\gamma_{\delta})/\delta$ outside the range $(0.5,2)$ were truncated to that range in Figure~\ref{fig:transition}.
As specified by our theory, the normalized exponent is $m=\log(1/\gamma_\delta)/\log(1/\delta)$ for $c>1/2$ and $m=\log(1/(1-\gamma_\delta))/\log(1/\delta)$ for $c<1/2$; it falls in the range delimited by $m_*$ and $m^*$.

\subsection{Calibration Algorithms}\label{sim:np_binary}

We investigate the performance of the calibration algorithms for binary classification and hypothesis testing detailed in Sections~\ref{sec:binary_calib} and~\ref{sec:np_calib}. These algorithms calibrate classifiers that estimate the regression function $\eta$, empirically targeting selective error rates in binary classification and Type~I and Type~II errors in hypothesis testing, as implemented in Algorithms~\ref{Binary_Class_Calib:alg} and~\ref{NP_Indec_Calib:alg}.

We consider two simulation settings. Simulation~1 uses $1{,}000$ training, calibration, and test observations drawn from a balanced Gaussian mixture with common variance $1$. The separation parameter is $\Delta = |\mu_1 - \mu_2|/2$, where smaller values of $\Delta$ correspond to lower-signal problems and $\Delta = 0$ implies that the two classes are indistinguishable. Simulation~2 uses $1{,}000$ training, calibration, and test observations drawn from a two-dimensional diamond data-generating mechanism whose Bayes decision boundary is diamond-shaped. We index this setting by $\Delta'$, with larger values of $\Delta'$ producing sharper class separation around the nonlinear boundary. For each value of $\Delta$ or $\Delta'$, results are averaged over $1{,}000$ Monte Carlo repetitions.

\begin{figure}[!t]
    \centering

    \scalebox{0.7}{\input{Figures/sim1_and2_binary}}

    \caption{\label{fig:binary_class_sims}\small Binary classification comparing classification with reject (orange / square), selection with conformal $p$-values (blue / triangle), our approach (green / circle), and an oracle procedure (gray / diamond) that has knowledge of the true class labels. Top: Gaussian mixture with separation $\Delta$. Bottom: Harder setting where the Bayes decision boundary is a diamond-shaped nonlinear curve. The left panel in each row shows the selective error rate among all definitive decisions, with the dashed line marking the target level $\alpha = 0.10$. The right panel shows the corresponding proportion of indecisions. Methods are compared across increasing signal strength $\Delta$ or $\Delta'$, illustrating both error control and the cost in undecided observations.}
\end{figure}

For the binary-classification simulations, we compare our calibrated indecision method with classification with reject, selection using conformal $p$-values, and an oracle benchmark. The two competing methods are adapted to target the same selective error-rate objective as our method, although this is not their original purpose. Classification with reject is implemented using the reject option hinge-loss as described in \citep{bartlett08a} and its tuning parameter is selected on the calibration sample to achieve the ydesired empirical selective error rate while making as many definitive decisions as possible. Similarly, the conformal procedure was originally designed for outlier detection, as described in \citep{Jin_23}, rather than selective classification. To adapt it to the binary setting, we apply the conformal construction separately to each class at the nominal level $2\cdot\alpha$, effectively reflecting the procedure across the two class labels. This produces class-specific conformal $p$-values for whether an observation is compatible with each class. Observations accepted by only one class are assigned to that class, observations accepted by neither class are assigned to the indecision set, and observations accepted by both classes are assigned to the class with the stronger confidence score. This adaptation allows the conformal method to be evaluated under the same selective error-control criterion as Algorithm \ref{Binary_Class_Calib:alg}, even though this was not its original objective.

\subsubsection{Binary Classification}\label{sim:binary}

The top row of Figure~\ref{fig:binary_class_sims} evaluates Algorithm~\ref{Binary_Class_Calib:alg} for a Gaussian mixture at the target error level $\alpha = 10\%$. The left panel compares our method (green / circle), using a GAM base learner \citep{Hastie2017}, with the oracle rule (black / dashed), which uses the true regression function $\eta$ and the true test labels to select the indecision threshold. All four methods have mean selective error at or below the target once $\Delta$ is larger than approximately $0.5$. For larger separations, around $\Delta \geq 1.5$, the problem becomes sufficiently easy that the methods have conservative mean error with almost no indecisions. When $\Delta < 0.5$, the problem is substantially harder, and nearly all observations must be assigned to the indecision set to maintain the desired error rate. In this low-signal regime, the methods behave similarly. The observed discrepancies in selective error are largely driven by the very small number of observations that receive definitive decisions.

\begin{figure}[!t]
    \centering

    \scalebox{0.8}{\input{Figures/sim1_and2_np}}

    \caption{\label{fig:np_sims}\small Hypothesis testing comparing the original NP procedure with one that uses indecisions. The first three columns show overall error, Type~I error, and Type~II error; dashed lines mark the target error level $\alpha = 0.10$ for Type~I and Type~II errors. The green shaded regions show the range of achievable error rates under varying indecision levels, and the green line shows the selected indecision-calibrated rule. The final column reports the selected proportion of indecisions. All other aspects are identical to the simulations for binary classification.}
\end{figure}

The bottom row of Figure~\ref{fig:binary_class_sims} evaluates the same binary-classification calibration problem as Simulation~1, but under the more challenging diamond data-generating mechanism. Here the data are generated in two dimensions so that the Bayes decision boundary is diamond-shaped, with larger~$\Delta'$ producing sharper class separation around this nonlinear boundary. In this setting, nearly all methods have mean selective error at or below the target once $\Delta' > 3$. The classification-with-reject method has more difficulty, with its mean error remaining slightly above the target in parts of the grid, reflecting the nonlinear decision boundary. Selection with conformal $p$-values meets the target for approximately $\Delta' > 2$, whereas our method does so over a somewhat broader range, beginning around $\Delta' = 1.5$. This improvement is most visible in the harder low-signal regimes, where the constrained Bayes motivation behind our calibration procedure yields a more robust selective rule than the conventional alternatives.

\subsubsection{NP Testing Paradigm} \label{sim:np}
For hypothesis testing, Algorithm~\ref{NP_Indec_Calib:alg} targets both Type~I and Type~II errors. We begin with a standard Neyman--Pearson (NP) classifier whose decision threshold is selected using a held-out calibration data set. Among the thresholds for which the estimated Type~I error is at or below $\alpha_1$, we select the threshold having the smallest estimated Type~II error. We then apply our indecision calibration step to the fitted NP classifier to target both error types among observations receiving definitive decisions. Throughout, we target $\alpha_1 = \alpha_2 = 0.10$.

Figure \ref{fig:np_sims} reports the NP results for the same two simulation settings used in the binary-classification study. The top row corresponds to the Gaussian simulation, where the NP classifier is fitted using LDA \citep{Hastie2017}. The bottom row corresponds to the diamond simulation, where the NP classifier is fitted using logistic regression on the nonlinear feature $|X_1| + |X_2|$, which captures the diamond-shaped Bayes boundary. In each row, the orange curve shows the original NP classifier, while the green curve shows our indecision-calibrated version of the same classifier. The first three panels report overall error, Type~I error, and Type~II error, respectively, and the rightmost panel shows the selected indecision mass $\gamma$.

In both simulations, the original NP classifier empirically meets the Type~I target and is selected to have the smallest estimated Type~II error among the candidate thresholds satisfying this constraint on the calibration data. However, it does not directly enforce the Type~II error target. Our method adds an indecision region to the NP classifier and selects $\gamma$ so that both error types are empirically at or below their targets. The green shaded regions show the range of achievable error rates over different values of $\gamma$, while the green curve corresponds to the selected value. In the harder low-signal regimes, larger indecision rates are needed to meet both targets; as the signal strength increases, the selected $\gamma$ decreases rapidly, indicating that fewer observations need to be left undecided once the testing problem becomes easier.

Finally, in these simulation settings, the use of indecisions almost always reduces the classifier's Type~II error relative to the baseline NP classifier. This observation should not be interpreted as a guarantee for arbitrary fitted rules or applications.

\section{Real Data}\label{sim:real_data_compas}

\begin{figure}[!t]
    \centering
    \scalebox{0.8}{\input{Figures/compas_binary_np}}

    \caption{\label{fig:compass_combined}\small Predicting criminal recidivism on the COMPAS dataset using binary selective classification and NP hypothesis testing with indecisions. \textbf{Top:} Binary classification results. The left panel shows the error rate of our data-driven method across fixed indecision masses $\gamma$; the middle panel compares classification with reject (orange / square), selection with conformal $p$-values (blue / triangle), and our data-driven method (green / circle) across target error levels $\alpha$; and the right panel shows the corresponding selected indecision proportions. \textbf{Bottom:} NP results comparing the unadjusted NP algorithm that cannot use indecisions with our approach that can utalize indecisions for hypothesis testing. The panels show overall error, Type~I error, Type~II error, and the selected indecision mass $\gamma$. Green shading indicates the 90\% range of achievable error rates under varying indecision levels.
    }
\end{figure}

Predicting criminal recidivism is a well-studied application of automated decision systems. A prominent example is the COMPAS risk-assessment instrument, which is used in the United States to assess a defendant's likelihood of reoffending. Given the high stakes of this task, the reliability of definitive predictions is critical, particularly when the underlying classification problem is difficult.

We analyze a dataset originally collected by ProPublica to investigate fairness in machine learning \citep{Angwin_16}. Although fairness is not our focus here, our methodology could be extended to provide group-specific error control for known protected groups. We approach the task using both binary classification and hypothesis testing, as outlined in Section \ref{sec:gen_hyp_test}, applying our data-driven Algorithms~\ref{Binary_Class_Calib:alg} and~\ref{NP_Indec_Calib:alg}.

The processed dataset contains $6,172$ defendants, of whom $2,990$ are labeled as recidivating under the outcome used in our analysis. We perform $100$ random splits of the data, assigning $25\%$ of the observations to training, $25\%$ to calibration, and the remaining $50\%$ to a held-out test set. For each split, we fit an ensemble consisting of AdaBoost, logistic regression, a generalized additive model (GAM), and nonparametric naive Bayes. Each component estimates the conditional probability of class~1, corresponding to no recorded recidivism, and the four estimated probabilities are averaged with equal weights. The resulting ensemble score is used for both the binary-classification and NP analyses. All classification thresholds and indecision regions are selected using only the calibration subset, after which error rates and indecision masses are evaluated on the held-out test subset and averaged across the $100$ splits \citep{Hastie2017}.

Figure~\ref{fig:compass_combined} summarizes the COMPAS real-data analysis. The top row reports the binary selective-classification results. Both our data-driven method and the conformal $p$-value method use scores obtained from the four-model ensemble. The left panel varies the indecision mass $\gamma$ for our data-driven method, illustrating the tradeoff between assigning more observations to the indecision set and reducing the error rate among definitive decisions. The middle panel varies the target error level $\alpha$ and compares classification with reject, selection with conformal $p$-values, and our data-driven method. The three procedures track the target more closely at moderate values of $\alpha$. At the most stringent target levels, however, the held-out error rates exceed their nominal targets because the procedures make definitive predictions for only a small proportion of observations. Our method nevertheless has the lowest held-out selective error among the three procedures at $\alpha=0.05$. The right panel shows the corresponding proportion of indecisions. As $\alpha$ increases, the error constraint becomes less stringent and fewer observations are assigned to the indecision set.

The bottom row reports the NP hypothesis-testing results. Both the calibration-selected NP classifier and our indecision-calibrated NP method use the same four-model ensemble score. For the baseline NP classifier, the decision threshold is selected on the calibration subset. Among the thresholds whose estimated Type~I error is at or below the target level, we select the threshold having the smallest estimated Type~II error. The ensemble is fitted using the COMPAS covariates age, sex, number of prior offenses, days between screening and arrest, and COMPAS decile score. The first three panels report overall error, Type~I error, and Type~II error, respectively, while the final panel reports the selected indecision mass $\gamma$. The orange curve represents the calibration-selected NP classifier, and the green curve represents our indecision-calibrated NP method. The green shading indicates the 5th--95th percentile range of achievable held-out error rates across the candidate indecision levels.

The baseline NP classifier is calibrated to keep its estimated Type~I error at or below the target while minimizing its estimated Type~II error subject to this constraint, but it does not directly impose a user-specified Type~II error target. Introducing indecisions substantially reduces Type~II error relative to this baseline. For moderate target levels, the held-out Type~I and Type~II errors of our method closely track their nominal values. At the most stringent levels, however, the selected rule makes very few definitive decisions, and the empirical calibration constraints do not translate into exact held-out control. For example, at the nominal level $\alpha_1=\alpha_2=0.05$, the average held-out Type~I and Type~II errors are approximately $0.096$ and $0.121$, respectively, while approximately $98\%$ of observations are assigned to the indecision set. As the target error level increases, the selected indecision mass decreases, matching the qualitative pattern observed in the simulations.

\bibliographystyle{chicago}
\bibliography{myrefs.bib}

\newpage

\setcounter{page}{1}

\appendix 

\renewcommand{\theequation}{\thesection.\arabic{equation}}
\setcounter{equation}{0}

\begin{center}
\Large Ask for More Than Bayes Optimal: A Theory of Indecisions for Selective Hypothesis Testing \\ Supplementary Material 
\end{center}
\begin{center}
    By: Mohamed Ndaoud$^1$, Peter Radchenko$^2$, and Bradley Rava$^2$
\end{center}

\footnotetext[1]{Department of Decision Sciences, ESSEC Business School, \url{ndaoud@essec.edu}}
\footnotetext[2]{Discipline of Business Analytics, University of Sydney Business School, \url{peter.radchenko@sydney.edu.au}, \url{bradley.rava@sydney.edu.au}}

\medskip
\medskip
\medskip

\section{Main Proofs}
\begin{lemma}\label{lem:constrained:min}
    Let $f$ and $g$ be two positive functions and $c>0$ and let $\mathcal{H}_c = \{A: \int_{A} f = c \}$. Assuming that $\mathcal{H}_c$ is not empty and that any required boundary mass below is attainable by a measurable subset, we have that any
    $$
   A_c^* \in \underset{A \in \mathcal{H}_c}{\arg\min} \int_A  g
    $$
    is of the form  $A_c^* := \{x:g(x)<t_c f(x)\}\cup\mathcal{M}_c$ for some $t_c \geq 0$ where $\mathcal{M}_c\subset\{x:g(x)=t_c f(x)\}$ such that $\int_{A_c^*} f = c$.
    
    In particular, if for all $t$, $\mu\{x:g(x)=t f(x)\}=0$, then $A_c^*$ is unique up to $\mu$-null sets and $A_c^* = \{x:g(x)\leq t_c f(x)\}$ almost surely.
\end{lemma}
\begin{proof}
    Observe that we may assume that $f>0$, since for any $A \in \mathcal{H}_c$ we also have that $B=A\cap\{x:f(x)>0\}\in\mathcal{H}_c$ and $\int_A  g \geq \int_B  g $. For the sake of generality we consider $f$ and $g$ to be simply positive.
    
    First, suppose that $A_c^*$ has the displayed threshold form. Then, for any $A \in \mathcal{H}_c$, we have
    \begin{align*}
        \int_A g-\int_{A^*_c}g &= \int_{A\setminus A^*_c}g-\int_{A^*_c\setminus A}g \\
        &\geq t_c\int_{A\setminus A^*_c}f-t_c\int_{A^*_c\setminus} Af \\
        &=t_c\left(\int_A f-\int_{A^*_c}f\right)=0.
    \end{align*}
    It follows that
    $$
   A_c^* \in \underset{A \in \mathcal{H}_c}{\arg\min} \int_A  g.
    $$
    We next show that, for any $c$, there exists $A_c^* := \{x:g(x)<t_c f(x)\}\cup\mathcal{M}_c$ for some $t_c \geq 0$ where $\mathcal{M}_c\subset\{x:g(x)=t_c f(x)\}$ and such that $\int_{A_c^*} f = c$. Define the nondecreasing function
    $$
    h(t):=\int_{\{x:g(x)\leq t f(x)\}} f,\qquad t\geq0.
    $$
    Hence, we can define, for any $m\geq0$,
    $$
    h^{-1}(m) = \inf\{t:h(t)\geq m\}.
    $$
    Let us set $t_c = h^{-1}(c)$ and $\mathcal{F}_c:=\{x:g(x)=t_c f(x)\}$.

    If $h(t_c) = c$, then we are done with $\mathcal{M}_c = \mathcal{F}_c$. Otherwise $h(t_c)>c$ and for any $t<t_c$, $h(t)<c$. In particular, $c^{-}:=\lim_{t\to t_c^{-}}h(t)\leq c$ and $h(t_c) - c^{-} = \int_{\mathcal{F}_c} f$.
    
    By the boundary-mass attainability assumption, there exists $\mathcal{M}_c\subseteq\mathcal{F}_c$ such that
    $$c - c^{-} = \int_{\mathcal{M}_c} f.$$
    By setting $A^*_c=\{x:g(x)<t_c f(x)\}\cup\mathcal{M}_c$, it follows that
    $$
    \int_{A_c^*} f = c^{-} + \int_{\mathcal{M}_c} f = c.
    $$

    It remains to show that any minimizer $A^*$ satisfies almost surely
    $$
    \{x:g(x)<t_c f(x)\}\subset A^*\subset\{x:g(x)\leq t_c f(x)\}.
    $$
Let us use the following notation $B_1=\{x:g(x)>t_c f(x)\}$, $B_2=\{x:g(x)<t_c f(x)\}$, and $B_3=\{x:g(x)=t_c f(x)\}$. In that case, we have
$$
\int_{A^*}  g = \int_{A^*\cap B_1}  g + \int_{A^*\cap B_2}  g + \int_{A^*\cap B_3}  g.
$$
It follows that
$$
0=\int_{A^*}g-\int_{A_c^*}g=\int_{A^*\cap B_1}g+\int_{(A^*\setminus A_c^*)\cap B_3}g-\int_{B_2\setminus} A^*g-\int_{(A_c^*\setminus A^*)\cap B_3}g.
$$
Similarly we also have that 
    $$
0=t_c\left(\int_{A^*}f-\int_{A_c^*}f\right)=t_c\int_{A^*\cap B_1}f+t_c\int_{(A^*\setminus A_c^*)\cap B_3}f-t_c\int_{B_2\setminus A^*}f-t_c\int_{(A_c^*\setminus A^*)\cap B_3}f.
$$
Combining both equations and the fact that $g(x) = t_c \cdot f(x)$ on $B_3$ leads to 
$$
0=\int_{A^*\cap B_1}(g-t_c f)-\int_{B_2\setminus A^*}(g-t_c f).
$$
Both terms on the right are nonnegative; hence both vanish. Therefore, $A^*\cap B_1$ and $B_2\setminus A^*$ are $\mu$-null, which proves the claimed inclusions.
\end{proof}

\subsection{Proof of Theorem \ref{thm:acc:lower_bound}}
    Let us consider a classifier $\tilde Y(X)$ such that $\mathbf{P}(\tilde Y = 0) = \gamma$ and let $A$ be the set where $\tilde{Y} = 0$. We have that
    \begin{align*}
    {\mathbf P}_Y(\tilde Y \neq Y | \tilde Y \neq 0 ) &= \frac{p_1 \mathbf{P}_1(\tilde Y = 2) + p_{2} \mathbf{P}_{2}(\tilde Y = 1)}{1 - \mathbf{P}(\tilde Y = 0)}\\
    & = \frac{\int_{A^c}\!\left[\mathbf{1}\{\tilde Y(x)=2\}p_1f_1(x)+\mathbf{1}\{\tilde Y(x)=1\}p_2f_2(x)\right]d\mu(x)}{1-\gamma}.
    \end{align*}
For each $x$, the integrand is minimized by choosing class~1 when $p_1f_1(x)\geq p_2f_2(x)$ and class~2 otherwise. Hence
$$
    {\mathbf P}_Y(\tilde Y \neq Y | \tilde Y \neq 0 ) \geq  \frac{\int_{A^c}(p_1f_1\wedge p_2f_2)\,d\mu}{1-\gamma}.
$$
Invoking Lemma \ref{lem:constrained:min} we get further that the above quantity is minimized for 
$$
A^* := \left\{x \; ; \; \frac{ p_1 f_1 \wedge p_{2} f_{2} }{p_1 f_1 + p_{2} f_{2}}(x) > 1 - \tau_\gamma  \right\} \cup \mathcal{M}_\gamma,
$$
such that $\mathbf{P}(A^*) = \gamma$ and $\tau_\gamma \in [1/2,1]$. The result and the expression for $Y^*$ follow.
\subsection{Proof of Proposition \ref{prop:risk}}
    The first part is straightforward from Theorem \ref{thm:acc:lower_bound}. Next observe that $\tau_\gamma$ is increasing by definition. Moreover $\underset{\gamma \to 1^{-}}{\lim } \tau_\gamma \leq 1$ and $\underset{\gamma \to 0^{+}}{\lim } \tau_\gamma \geq  1/2$.

    On the one hand, let $\beta_1 \geq \beta_2$ and hence $\tau_{\beta_1} \geq \tau_{\beta_2}$. We have
    \begin{align*}
        \mathcal{R}(\beta_2) &= \frac{\int_{A_{\beta_2}^{*c}}  (p_1 f_1 \wedge p_{2} f_{2})}{1-\beta_2} \\
        &= \frac{\int_{A_{\beta_1}^{*c}}  (p_1 f_1 \wedge p_{2} f_{2})}{1-\beta_2} + \frac{\int_{A_{\beta_2}^{*c}\setminus A_{\beta_1}^{*c}}  (p_1 f_1 \wedge p_{2} f_{2})}{1-\beta_2} \\
        &\geq \frac{1}{1-\beta_2} \left( \int_{A_{\beta_1}^{*c}}  (p_1 f_1 \wedge p_{2} f_{2}) + \frac{\int_{A_{\beta_1}^{*c}}  (p_1 f_1 \wedge p_{2} f_{2})}{1 - \beta_1}  \int_{A_{\beta_2}^{*c}\setminus A_{\beta_1}^{*c}}  (p_1 f_1 + p_{2} f_{2}) \right),
    \end{align*}
    where we have used, in the last inequality, the fact that on $A_{\beta_2}^{*c}\setminus A_{\beta_1}^{*c}$ we have that 
    $$
      \frac{(p_1 f_1 \wedge p_{2} f_{2})}{(p_1 f_1 + p_{2} f_{2})} \geq (1-\tau_{\beta_1}),
    $$
    while on $A_{\beta_1}^{*c}$ we have
    $$
     \frac{(p_1 f_1 \wedge p_{2} f_{2})}{(p_1 f_1 + p_{2} f_{2})} \leq (1-\tau_{\beta_1}).
    $$
    As a consequence we have that 
$$
\frac{(p_1 f_1 \wedge p_{2} f_{2})}{(p_1 f_1 + p_{2} f_{2})} \geq \frac{\int_{A_{\beta_1}^{*c}}  (p_1 f_1 \wedge p_{2} f_{2})}{1 - \beta_1}
$$
    on $A_{\beta_2}^{*c}\setminus A_{\beta_1}^{*c}$. It follows that
    $$
    \mathcal{R}(\beta_2) \geq \frac{\int_{A_{\beta_1}^{*c}}  (p_1 f_1 \wedge p_{2} f_{2})}{1-\beta_2} \frac{1-\beta_2}{1-\beta_1} \geq \mathcal{R}(\beta_1).
    $$
    Consequently, $\mathcal{R}(\gamma)$ is non-increasing.
    
    On the other hand, we have that
    $$
    \mathcal{R}(\beta_2) - \mathcal{R}(\beta_1) =  \frac{(\beta_2 - \beta_1)\int_{A_{\beta_1}^{*c}}  (p_1 f_1 \wedge p_{2} f_{2})}{(1-\beta_2)(1-\beta_1)} + \frac{\int_{A_{\beta_2}^{*c}\setminus A_{\beta_1}^{*c}}  (p_1 f_1 \wedge p_{2} f_{2})}{1-\beta_2}.
    $$
    On the event $A_{\beta_2}^{*c}\setminus A_{\beta_1}^{*c}$ we have that
    $$
    p_1 f_1 \wedge p_{2} f_{2} \leq (p_1 f_1 + p_{2} f_{2})(1-\tau_{\beta_2}).
    $$
    Hence
    $$
\int_{A_{\beta_2}^{*c}\setminus A_{\beta_1}^{*c}}  (p_1 f_1 \wedge p_{2} f_{2}) \leq (1-\tau_{\beta_2}) (\beta_1 - \beta_2).
    $$
    As a consequence we get that
    $$
    0 \leq \mathcal{R}(\beta_2) - \mathcal{R}(\beta_1) \leq \frac{2(\beta_1 - \beta_2)}{1-\beta_2}.
    $$
    
    We conclude that $\gamma \to \mathcal{R}(\gamma)$ is continuous. This proof is complete.

\subsection{Proof of Theorem \ref{thm:NP:lower_bound}}

We proceed in two steps. First, let us consider a reparametrization of the problem through the function $\Gamma$, defined as
$$
\Gamma(\gamma,r) := \inf_{\tilde Y_{\gamma}} \frac{{\mathbf P}_{2}(\tilde Y_{\gamma} = 1)}{\mathbf P_2(\tilde Y_{\gamma} \neq 0)},
$$
where the infimum is taken over all classifiers $\tilde{Y}_{\gamma}$ taking values in $\{0,1,2\}$ such that ${\mathbf P}(\tilde Y_{\gamma} = 0) = \gamma$ and ${\mathbf P}(\tilde Y_{\gamma} = 2) = r$. It is straightforward to see that this minimization problem is equivalent to minimizing the ratio $\frac{{\mathbf P}_{2}(\tilde Y_{\gamma} = 1)}{\mathbf P_2(\tilde Y_\gamma = 2)}$. Given the two constraints ${\mathbf P}(\tilde Y_\gamma = 2) = r$ and ${\mathbf P}(\tilde Y_\gamma = 1) = 1-\gamma- r$, we can decompose this problem into two separate optimization tasks: the first minimizes ${\mathbf P}_{2}(\tilde Y_{\gamma} = 1)$ subject to ${\mathbf P}(\tilde Y_\gamma = 1) = 1-\gamma- r$, and the second maximizes $\mathbf P_2(\tilde Y_\gamma = 2)$ subject to ${\mathbf P}(\tilde Y_\gamma = 2) = r$. 

Invoking Lemma \ref{lem:constrained:min}, these quantities are optimized by choosing the decision regions
\[
A^*_1=\{\eta>\tau_2\}\cup \mathcal M^2_{\gamma,r},
\qquad
A^*_2=\{\eta\le\tau_1\}\setminus \mathcal M^1_{\gamma,r},
\]
such that $\mathbf{P}(A_1^*) = 1-\gamma-r$ and $\mathbf{P}(A_2^*) = r$. The explicit expression for the optimal classifier $Y_\gamma^*$ follows immediately, yielding
$$
\Gamma(\gamma,r) = \frac{\int_{A^*_1} f_{2}\,d\mu} {\int_{A^*_2} f_{2}\,d\mu + \int_{A^*_1} f_{2}\,d\mu}.
$$

In the second step, we map the dependence on $r$ to a dependence on the Type~I error, denoted by $\alpha$. Analogously to standard classification settings, $\tau_2$ is  increasing with $\gamma + r$, while $\tau_1$ is  increasing with $r$. Note that $\tau_2 = H^{-1}(r+\gamma)$ and $\tau_1 = H^{-1}(r)$, where $H(t):=\mathbf P\{\eta(X)\leq t\}$ is the CDF of the posterior score under the mixture distribution.

Recall that the Type~I error corresponding to the optimal solution for a given pair $(\gamma,r)$ is given by
$$
\mathcal{P}^r_{\text{I}} := \frac{\int_{A^*_2} f_{1}\,d\mu} {\int_{A^*_2} f_{1}\,d\mu + \int_{A^*_1} f_{1}\,d\mu}.
$$
For a fixed $\gamma$, the numerator is nondecreasing with $r$, while the second term of the denominator ($\int_{A^*_1} f_{1}\,d\mu$) is  decreasing with $r$. Consequently, $\mathcal{P}^r_{\text{I}}$ is nondecreasing with $r$.

To establish continuity, let $w_1(r):= \int_{A^*_2} f_{1}\,d\mu$ and $w_2(r):=\int_{A^*_1} f_{1}\,d\mu$. Recall that $p_1\int_{A^*_1} f_{1}\,d\mu + p_2\int_{A^*_1} f_{2}\,d\mu = 1-\gamma-r$ and $p_1\int_{A^*_2} f_{1}\,d\mu + p_2\int_{A^*_2} f_{2}\,d\mu = r$. Choosing $r < r'$, the tail constraints imply
$$
0 \leq w_1(r') - w_1(r) \leq \frac{r'-r}{p_1},
$$
and similarly,
$$
0 \leq w_2(r) - w_2(r') \leq \frac{r'-r}{p_1}.
$$
Since both mappings are Lipschitz continuous, $\mathcal{P}^r_{\text{I}}$ is continuous and nondecreasing. Thus, for any admissible $\alpha_1$, there exists an $r_{\alpha_1} \in [0,1-\gamma]$ such that
$$
\alpha_1 = \mathcal{P}^{r_{\alpha_1}}_{\text{I}}.
$$
We choose $r_{\alpha_1}$ as the supremum of the values of $r$ satisfying $\mathcal{P}^{r}_{\text{I}}\leq\alpha_1$.

Finally, we claim that 
$$
\mathcal{P}(\alpha_1,\gamma) = \Gamma(\gamma,r_{\alpha_1}).
$$
Suppose for contradiction that $\mathcal{P}(\alpha_1,\gamma) < \Gamma(\gamma,r_{\alpha_1})$. Then, there exists an optimal set $A_2^{**}$ such that the Type~II error is strictly less than $\Gamma(\gamma,r_{\alpha_1})$, the Type~I error is less than $\alpha_1$ and where $\mathbf{P}(A_2^{**}):= r' > r_{\alpha_1}$. This implies
$$
\mathcal{P}^{r'}_{\text{I}} \leq \alpha_1.
$$
Since $r_{\alpha_1}$  is the largest of all $r$ such that $ \alpha_1 \geq  \mathcal{P}^{r}_{\text{I}}$, this necessitates $r_{\alpha_1} \geq r'$, a clear contradiction.

\subsection{Proof of Proposition \ref{prop:NP}}

We now allow the parameter $\gamma$ to vary. Consequently, the choice of $r_\alpha$ depends on $\gamma$, satisfying $\alpha_1 = \mathcal{P}^{r_{\alpha_1}(\gamma)}_{\text{I}}$. Since $\alpha_1$ is fixed throughout this section, we simplify notation by writing $r_\gamma := r_{\alpha_1}(\gamma)$.

Recall that the optimal decision regions are defined as
\[
A^*_1=\{\eta>\tau_2\}\cup \mathcal M^2_{\gamma,r_\gamma},
\qquad
A^*_2=\{\eta\le\tau_1\}\setminus \mathcal M^1_{\gamma,r_\gamma},
\]
such that $\mathbf{P}(A_1^*) = 1-\gamma-r_\gamma$ and $\mathbf{P}(A_2^*) = r_\gamma$. Under this setup, we have
$$
\Gamma(\gamma,r_\gamma) = \frac{\int_{A^*_1} f_{2}\,d\mu} {\int_{A^*_2} f_{2}\,d\mu + \int_{A^*_1} f_{2}\,d\mu},
$$
and
$$
\alpha_1 = \frac{\int_{A^*_2} f_{1}\,d\mu} {\int_{A^*_2} f_{1}\,d\mu + \int_{A^*_1} f_{1}\,d\mu}.
$$
Moreover, the thresholds satisfy $\tau_2 = H^{-1}(r_\gamma +\gamma)$ and $\tau_1 = H^{-1}(r_\gamma)$.

We first observe that the mapping $\gamma \mapsto r_\gamma$ is  decreasing, whereas $\gamma \mapsto \gamma + r_\gamma$ is  increasing. Indeed, since $\alpha_1$ is fixed, the ratio $\int_{A^*_2} f_{1}\,d\mu / \int_{A^*_1} f_{1}\,d\mu$ must remain constant as $\gamma$ increases. This invariance requires $\tau_1$ and $\tau_2$ to exhibit opposite monotonic behaviors. If $r_\gamma$ were to increase, both $\tau_1$ and $\tau_2$ would simultaneously increase, which is impossible. Thus, $r_\gamma$ must be  decreasing. This implies that $\tau_1$ is decreasing, which forces $\tau_2$ to be increasing, meaning that $\gamma + r_\gamma$ is  increasing.

Next, notice that $\Gamma(\gamma,r_\gamma)$ shares the same monotonicity (with respect to $\gamma$) as the ratio $\int_{A^*_1} f_{2}\,d\mu / \int_{A^*_2} f_{2}\,d\mu$. Because the Type~I error is held constant, the monotonicity of $\Gamma(\gamma,r_\gamma)$ matches that of the product:
$$
\frac{\int_{A^*_1} f_{2}\,d\mu}{\int_{A^*_2} f_{2}\,d\mu} \frac{\int_{A^*_2} f_{1}\,d\mu}{\int_{A^*_1} f_{1}\,d\mu} = \mathbf{E}_1\left[\frac{f_2}{f_1}(X) \,\middle|\, X \in A^*_1\right] \mathbf{E}_2\left[\frac{f_1}{f_2}(X) \,\middle|\, X \in A^*_2\right].
$$

Observing that $f_2/f_1 = (1/\eta - 1 ) p_1/p_2$, the likelihood ratio $f_2/f_1$ is strictly decreasing with $\eta$. As $\gamma$ increases, $\tau_2$ increases, causing the set $A^*_1$ to shrink monotonically (in the set-inclusion sense). Consequently, the conditional expectation $\mathbf{E}_1\left[f_2/f_1(X) \mid X \in A^*_1\right]$ decreases. By an identical argument, $\mathbf{E}_2\left[f_1/f_2(X) \mid X \in A^*_2\right]$ also decreases as $\gamma$ increases. We conclude that $\gamma \mapsto \Gamma(\gamma,r_\gamma)$ is  decreasing.

For continuity, let $w_3(\gamma):= \int_{A^*_1} f_{2}\,d\mu$ and $w_4(\gamma):=\int_{A^*_2} f_{2}\,d\mu$. Recall that $p_1\int_{A^*_1} f_{1}\,d\mu + p_2\int_{A^*_1} f_{2}\,d\mu = 1-\gamma-r_\gamma$ and $p_1\int_{A^*_2} f_{1}\,d\mu + p_2\int_{A^*_2} f_{2}\,d\mu = r_\gamma$. For $\gamma < \gamma'$, the tail constraints dictate
$$
0 \leq w_3(\gamma) - w_3(\gamma') \leq \frac{\gamma'-\gamma}{p_2},
$$
and
$$
0 \leq w_4(\gamma) - w_4(\gamma') \leq \frac{\gamma'-\gamma}{p_2}.
$$
Being Lipschitz continuous, these mappings are continuous, confirming that $\gamma \mapsto \Gamma(\gamma,r_\gamma)$ is continuous and nonincreasing. Consequently, for any achievable level $\alpha_2$, there exists at least one $\gamma\in[0,1)$ such that
$$
\alpha_2 = \mathcal{P}(\alpha_1,\gamma).
$$

\subsection{Proofs for Section~\ref{sec.plugin}}

\subsubsection{Proof of Lemma~\ref{lem.risk.expr}}

Note that
\begin{align}
(1-\gamma)\big[R(\widehat{Y}_{\gamma}) - R(Y^*_{\gamma})\big] 
&= E\eta\big(\mathbf{1}\{\widehat{Y}_{\gamma}=2, Y^*_{\gamma}\ne2\} - \mathbf{1}\{Y^*_{\gamma}=2, \widehat{Y}_{\gamma}\ne2\}\big) \nonumber \\
&\quad + E(1-\eta)\big(\mathbf{1}\{\widehat{Y}_{\gamma}=1, Y^*_{\gamma}\ne1\} - \mathbf{1}\{Y^*_{\gamma}=1, \widehat{Y}_{\gamma}\ne1\}\big). \label{lem1.prf.eq1}
\end{align}
Also note that 
\begin{eqnarray*}
&&\mathbf{1}\{\widehat{Y}_{\gamma}=2,Y^*_{\gamma}\ne2\}-\mathbf{1}\{Y^*_{\gamma}=2,\widehat{Y}_{\gamma}\ne2\}+\mathbf{1}\{\widehat{Y}_{\gamma}=1,Y^*_{\gamma}\ne1\}-\mathbf{1}\{Y^*_{\gamma}=1,\widehat{Y}_{\gamma}\ne1\}\\
&&~~~~~~~~~~~~~~~~~~~~~~~~~~~~~~~~~~~~~~~~~~~~~~~~~~~~~~~~~~~~~~~~~~~~~~~~~~~~~=\,\mathbf{1}\{Y^*_{\gamma}=0\}-\mathbf{1}\{\widehat{Y}_{\gamma}=0\}.
\end{eqnarray*}
Consequently, equality $\mathbf{P}(Y^*_{\gamma}=0)=\mathbf{P}(\widehat{Y}_{\gamma}=0)$ yields
\begin{align}
& E\big(\mathbf{1}\{\widehat{Y}_{\gamma}=2, Y^*_{\gamma}\ne2\} - \mathbf{1}\{Y^*_{\gamma}=2, \widehat{Y}_{\gamma}\ne2\}\big) \nonumber \\
&\quad + E\big(\mathbf{1}\{\widehat{Y}_{\gamma}=1, Y^*_{\gamma}\ne1\} - \mathbf{1}\{Y^*_{\gamma}=1, \widehat{Y}_{\gamma}\ne1\}\big) = 0. \label{lem1.prf.eq2}
\end{align}
Combining equations~(\ref{lem1.prf.eq1}) and~(\ref{lem1.prf.eq2}), we derive
\begin{eqnarray*}
(1-\gamma)\big[R(\widehat{Y}_{\gamma})-R(Y^*_{\gamma})\big]&=&E\big(\eta-[1-\tau^*_{\gamma}]\big)\big(\mathbf{1}\{\widehat{Y}_{\gamma}=2,Y^*_{\gamma}\ne2\}-\mathbf{1}\{Y^*_{\gamma}=2,\widehat{Y}_{\gamma}\ne2\}\big)\\
&~+&E\big(1-\eta-[1-\tau^*_{\gamma}]\big)\big(\mathbf{1}\{\widehat{Y}_{\gamma}=1,Y^*_{\gamma}\ne1\}-\mathbf{1}\{Y^*_{\gamma}=1,\widehat{Y}_{\gamma}\ne1\}\big).
\end{eqnarray*}
Note that $\eta<1-\tau^*_{\gamma}$ if and only if $Y^*_{\gamma}=2$. Also note that $\eta>\tau^*_{\gamma}$ if and only if $Y^*_{\gamma}=1$. Hence, we can rewrite the above display as follows:
\begin{eqnarray*}
(1-\gamma)\big[R(\widehat{Y}_{\gamma})-R(Y^*_{\gamma})\big]&=&E|1-\tau^*_{\gamma}-\eta|\big(\mathbf{1}\{\widehat{Y}_{\gamma}=2,Y^*_{\gamma}\ne2\}+\mathbf{1}\{Y^*_{\gamma}=2,\widehat{Y}_{\gamma}\ne2\}\big)\\
&~+&E|\tau^*_{\gamma}-\eta|\big(\mathbf{1}\{\widehat{Y}_{\gamma}=1,Y^*_{\gamma}\ne1\}+\mathbf{1}\{Y^*_{\gamma}=1,\widehat{Y}_{\gamma}\ne1\}\big),
\end{eqnarray*}
which completes the proof. \qed

\subsubsection{Proof of Theorem~\ref{lem.fixed.indec}} Define event~$A_{\phi}$ as follows:
\begin{eqnarray*}
    A_{\phi}&=&\{\tau^*_{\gamma}< \eta\le\widehat{\tau}_{\gamma}-\phi\}\cup\{\tau^*_{\gamma}< 1-\eta\le\widehat{\tau}_{\gamma}-\phi\}\\
    &&~~~~~~~~~~~~~~~~~~~~~~~~~~~~~~\cup\{\widehat{\tau}_{\gamma}+\phi< \eta\le \tau^*_{\gamma}\}\cup\{\widehat{\tau}_{\gamma}+\phi< 1-\eta\le \tau^*_{\gamma}\}.
\end{eqnarray*}
We will use the following result, which is proved in Section~\ref{sec:prf.lem.two-bounds}.
\begin{lemma}
\label{lem.two-bounds}
For $0<\phi\le\phi^*_{\gamma}$, we have
\begin{equation}
    \mathbf{P}\big(A_{\phi}\big)\le \mathbf{P}\big(|\widehat{\eta}-\eta|>\phi\big)\qquad\text{and}\qquad
    \mathbf{P}\big(|\widehat{\tau}_{\gamma}-\tau^*_{\gamma}|>2\phi\big)\lesssim\frac{\mathbf{P}\big(|\widehat{\eta}-\eta|>\phi\big)}{\phi^{\beta'}}.
\end{equation}
\end{lemma}
It follows from the proof of Lemma~\ref{lem.risk.expr}, that the equality in the statement of Lemma~\ref{lem.risk.expr} continues to hold when~$\tau^*_{\gamma}$ is replaced by an arbitrary constant~$c$. Moreover, another small modification to the proof allows us to replace~$c$  with~$\widehat{\tau}_{\gamma}$. We will focus on the first of the four terms in the resulting expression for $(1-\gamma)[R(\widehat{Y}_{\gamma})-R(Y^*_{\gamma})]$ -- the other three terms can be handled by analogous arguments. The term of interest can be bounded as follows:
$$
E|\widehat{\tau}_{\gamma}-\eta|\mathbf{1}\{Y^*_{\gamma}=1,\widehat{Y}_{\gamma}\ne Y^*_{\gamma}\}\le E_1+E_2+E_3,
$$
where
\begin{eqnarray*}
E_1&=&\mathbf{P}\big(\eta>\widehat{\tau}_{\gamma}+\phi,Y^*_{\gamma}=1,\widehat{Y}_{\gamma}\ne Y^*_{\gamma}\big),\\
E_2&=&\mathbf{P}\big(A_{\phi}\big),\quad\text{and}\\
E_3&=&\phi\mathbf{P}\big(|\eta-\widehat{\tau}_{\gamma}|\le\phi,Y^*_{\gamma}=1,\widehat{Y}_{\gamma}\ne Y^*_{\gamma}\big).
\end{eqnarray*}
Note that $\{Y^*_{\gamma}=1,\widehat{Y}_{\gamma}\ne Y^*_{\gamma}\}=
\{\eta>\tau^*_{\gamma},\widehat{\eta}\le\widehat{\tau}_{\gamma}\}$. Consequently, taking into account Lemma~\ref{lem.two-bounds}, we derive
\begin{equation}
    E_1+E_2\le 2\mathbf{P}\big(|\widehat{\eta}-\eta|>\phi\big). 
\end{equation}
We also have
\begin{align*}
E_3 &= \phi\, \mathbf{P}\big(|\widehat{\tau}_{\gamma} - \tau^*_{\gamma}| > 2\phi,\ |\eta - \widehat{\tau}_{\gamma}| \le \phi,\ \eta > \tau^*_{\gamma},\ \widehat{\eta} \le \widehat{\tau}_{\gamma}\big) \\
    &\quad + \phi\, \mathbf{P}\big(|\widehat{\tau}_{\gamma} - \tau^*_{\gamma}| \le 2\phi,\ |\eta - \widehat{\tau}_{\gamma}| \le \phi,\ \eta > \tau^*_{\gamma},\ \widehat{\eta} \le \widehat{\tau}_{\gamma}\big) \\
    &\le \phi\, \mathbf{P}\big(|\widehat{\tau}_{\gamma} - \tau^*_{\gamma}| > 2\phi\big)\, \mathbf{P}\big(\eta > \tau^*_{\gamma}\big) + \phi\, \mathbf{P}\big(|\eta - \tau^*_{\gamma}| \le 3\phi\big) \\
    &\lesssim \left(\phi^{1 - \beta'}\, \mathbf{P}\big(|\widehat{\eta} - \eta| > \phi\big) \wedge \phi\right)(1 - \gamma) + \phi^{1 + \beta},
\end{align*}
where we used Lemma~\ref{lem.two-bounds} and condition (\ref{lbnds_eta_prob01}) to derive the final bound. Thus, we get the desired bound for the first term in our resulting expression for $(1-\gamma)[R(\widehat{Y}_{\gamma})-R(Y^*_{\gamma})]$:
\begin{equation}
\label{fin.bnd1.prf.lem.fixed.indec}
    E|\widehat{\tau}_{\gamma}-\eta|\mathbf{1}\{Y^*_{\gamma}=1,\widehat{Y}_{\gamma}\ne Y^*_{\gamma}\}\lesssim \mathbf{P}\big(|\widehat{\eta}-\eta|>\phi\big)+(\phi^{1-\beta'}\mathbf{P}\big(|\widehat{\eta}-\eta|>\phi\big) \wedge \phi)(1-\gamma)+\phi^{1+\beta}.
\end{equation}
The other three terms can be similarly bounded using analogous arguments. This completes the proof of claim~(\ref{cl1.lem.fixed.indec}) in Theorem~\ref{lem.fixed.indec}. \qed
\subsubsection{Proof of Theorem~\ref{lem.fixed.acc}}
To simplify the notation, we will write ${Y^*}^{\tau}$ and~$\widehat{Y}^{\tau}$ for the classifiers~$Y^*$ and~$\widehat{Y}$, respectively, that use $\tau$ as the threshold. For example, 
$$
\widehat{Y}^{\tau}(X)=1\times\mathbf{1}\{ \widehat{\eta}(X)>{\tau}\}+2\times\mathbf{1}\{ \widehat{\eta}(X)<1-{\tau}\}.
$$
Recall that~$\gamma$ is the indecision level of the classifier $Y^*$ that uses threshold $\tau^*_{\gamma}$; also recall that $R^*=R({Y^*}^{\tau^*_{\gamma}})$.
We will use the following result, which is proved in Section~\ref{sec:prf.lem.bal}.
\begin{lemma}
\label{lem.bal}
For $\epsilon, \phi$ such that $0<\epsilon\le\phi\le \phi^*_{\gamma}\wedge(1/2-\tau^*_{\gamma}/2)$, we have
\begin{eqnarray*}
R^*-R({Y^*}^{\tau^*_{\gamma}+\phi})&\gtrsim& \phi^{1+2\beta'-\beta} \qquad \text{and}\\
R_{\widehat{\eta}}({\widehat{Y}}^{\tau^*_{\gamma}+\phi})-R({Y^*}^{\tau^*_{\gamma}+\phi})&\lesssim& \mathbf{P}\{ |\widehat{\eta}-\eta|>\epsilon\,|\,\widehat{\eta}\} +\epsilon\phi^{\beta}.
\end{eqnarray*}
\end{lemma}

\smallskip

\noindent By the definitions of~$\widehat{\tau}$ and~$\tau^*_{\gamma}$, the event $\{\widehat{\tau}>\tau^*_{\gamma}+\phi\}$ implies $\{R_{\widehat{\eta}}({\widehat{Y}}^{\tau^*_{\gamma}+\phi})>R^*\}$. Hence,
\begin{eqnarray*}
    \mathbf{P}\big( \widehat{\tau}>\tau^*_{\gamma}+\phi \big) &\le& 
    \mathbf{P}\Big( R_{\widehat{\eta}}({\widehat{Y}}^{\tau^*_{\gamma}+\phi})-R({Y^*}^{\tau^*_{\gamma}+\phi})>R^*-R({Y^*}^{\tau^*_{\gamma}+\phi})\Big).
\end{eqnarray*}

Using Lemma~\ref{lem.bal} to bound the components of the event on the right-hand side in the above display, we derive
\begin{eqnarray*}
    \mathbf{P}\big( \widehat{\tau}>\tau^*_{\gamma}+\phi \big) &\le& 
    \mathbf{P}\Big( \mathbf{P}\{ |\widehat{\eta}-\eta|>\epsilon\,|\,\widehat{\eta}\} +\epsilon\phi^{\beta}\gtrsim \phi^{1+2\beta'-\beta}\Big)\\
    &\le&
   \mathbf{P}\Big( \mathbf{P}\{ |\widehat{\eta}-\eta|>\epsilon\,|\,\widehat{\eta}\} \gtrsim \phi^{1+2\beta'-\beta}\Big) + 
    \mathbf{P}\big( \epsilon\phi^{\beta} \gtrsim \phi^{1+2\beta'-\beta}\big)\\
    &\lesssim&
   \frac{\mathbf{P}\{ |\widehat{\eta}-\eta|>\epsilon\}}{\phi^{1+2\beta'-\beta}} + 
    \mathbf{P}\big( \epsilon \gtrsim \phi^{1+2\beta'-2\beta}\big). 
\end{eqnarray*}
We take $\epsilon=c_1\phi^{1+2\beta'-2\beta}$ and note that we can choose $c_1$ sufficiently small to ensure that the second term in the line above is zero (recall that $\beta' \geq \beta$). This completes the proof of the first bound in Theorem~\ref{lem.fixed.acc}. The second bound in Theorem~\ref{lem.fixed.acc} follows from the first bound together with condition (\ref{lbnds_eta_prob01}). \qed

\subsubsection{Proof of Lemma~\ref{lem.two-bounds}}
\label{sec:prf.lem.two-bounds}
Note that~$\widehat{\tau}_{\gamma}$ is fully determined by~$\widehat{\eta}$. When $\widehat{\tau}_{\gamma}\ge\tau^*_{\gamma}+\phi$, we have
\begin{equation}
    \mathbf{P}\big(1-\widehat{\tau}_{\gamma}+\phi\le\eta\le\widehat{\tau}_{\gamma}-\phi\,|\,\widehat{\eta} \big)=\gamma+\mathbf{P}\big(\eta\in A_{\phi}\,|\,\widehat{\eta} \big).\label{Ad.bnd1}
\end{equation}
We also have 
\begin{eqnarray}
\mathbf{P}\big(1-\widehat{\tau}_{\gamma}+\phi\le\eta\le\widehat{\tau}_{\gamma}-\phi\,|\,\widehat{\eta} \big)&\le&\mathbf{P}\big(1-\widehat{\tau}_{\gamma}\le\widehat{\eta}\le\widehat{\tau}_{\gamma}\,|\,\widehat{\eta} \big)+\mathbf{P}\big(|\widehat{\eta}-\eta|>\phi\,|\,\widehat{\eta} \big)\nonumber\\
&=&\gamma+\mathbf{P}\big(|\widehat{\eta}-\eta|>\phi\,|\,\widehat{\eta} \big).\label{Ad.bnd2}
\end{eqnarray}
Combining~(\ref{Ad.bnd1}) and~(\ref{Ad.bnd2}), we derive
\begin{equation}
    \mathbf{P}\big(\eta\in A_{\phi}\,|\,\widehat{\eta} \big)\le\mathbf{P}\big(|\widehat{\eta}-\eta|>\phi\,|\,\widehat{\eta} \big).\label{Ad.bnd3}
\end{equation}

When $\tau^*_{\gamma}-\phi<\widehat{\tau}_{\gamma}<\tau^*_{\gamma}+\phi$, we have $\mathbf{P}\big(\eta\in A_{\phi}\,|\,\widehat{\eta} \big)=0$, and hence inequality~(\ref{Ad.bnd3}) still holds. Now consider the last remaining case:  $\widehat{\tau}_{\gamma}\le\tau^*_{\gamma}-\phi$. Note that 
\begin{equation}
    \mathbf{P}\big(1-\widehat{\tau}_{\gamma}-\phi\le\eta\le\widehat{\tau}_{\gamma}+\phi\,|\,\widehat{\eta} \big)=\gamma-\mathbf{P}\big(\eta\in A_{\phi}\,|\,\widehat{\eta} \big).\label{Ad.bnd4}
\end{equation}
We also have 
\begin{eqnarray}
\gamma&=&\mathbf{P}\big(1-\widehat{\tau}_{\gamma}\le\widehat{\eta}\le\widehat{\tau}_{\gamma}\,|\,\widehat{\eta} \big)\nonumber\\
&\le&\mathbf{P}\big(1-\widehat{\tau}_{\gamma}-\phi\le\eta\le\widehat{\tau}_{\gamma}+\phi\,|\,\widehat{\eta} \big)+\mathbf{P}\big(|\widehat{\eta}-\eta|>\phi\,|\,\widehat{\eta} \big).\label{Ad.bnd5}
\end{eqnarray}
Combining~(\ref{Ad.bnd4}) and~(\ref{Ad.bnd5}), we again derive inequality~(\ref{Ad.bnd3}), concluding that~(\ref{Ad.bnd3}) holds for all possible~$\widehat{\eta}$. Integrating~(\ref{Ad.bnd3}) over~$\widehat{\eta}$, we derive the first claim of Lemma~\ref{lem.two-bounds}.

For the second claim of Lemma~\ref{lem.two-bounds}, we will focus on bounding $\mathbf{P}\big(\widehat{\tau}_{\gamma}>\tau^*_{\gamma}+2\phi \big)$; the complementary bound on $\mathbf{P}\big(\widehat{\tau}_{\gamma}<\tau^*_{\gamma}-2\phi \big)$ follows analogously.  Note that $\widehat{\tau}_{\gamma}>\tau^*_{\gamma}+2\phi$ implies
\begin{eqnarray}
\mathbf{P}\big(1-\tau^*_{\gamma}-\phi \le\eta\le\tau^*_{\gamma}+\phi\big)&\le&
\mathbf{P}\big(1-\tau^*_{\gamma}-2\phi \le\widehat{\eta}\le\tau^*_{\gamma}+2\phi\,|\,\widehat{\eta}\big)+\mathbf{P}\big(|\widehat{\eta}-\eta|>\phi\,|\,\widehat{\eta} \big)\nonumber\\
&\le&
\gamma+\mathbf{P}\big(|\widehat{\eta}-\eta|>\phi\,|\,\widehat{\eta} \big).\label{tau-hat.bnd1}
\end{eqnarray}
By condition (\ref{lbnds_eta_prob01}), we also have
\begin{equation}
\label{tau-hat.bnd2}
 \mathbf{P}\big(1-\tau^*_{\gamma}-\phi \le\eta\le\tau^*_{\gamma}+\phi\big)\ge \gamma + c\phi^{\beta'},   
\end{equation}
for some fixed positive constant~$c$. Combining~(\ref{tau-hat.bnd1}) and~(\ref{tau-hat.bnd2}), we deduce that $\widehat{\tau}_{\gamma}>\tau^*_{\gamma}+2\phi$ implies $\mathbf{P}\big(|\widehat{\eta}-\eta|>\phi\,|\,\widehat{\eta} \big)\ge c\phi^{\beta'}$. Applying Markov's inequality, we then conclude that
\begin{equation*}
\mathbf{P}\big(\widehat{\tau}_{\gamma}>\tau^*_{\gamma}+2\phi\big)\le
\mathbf{P}\Big(\mathbf{P}\big(|\widehat{\eta}-\eta|>\phi\,|\,\widehat{\eta} \big)\ge c\phi^{\beta'}\Big)\le \frac{\mathbf{P}\big(|\widehat{\eta}-\eta|>\phi\big)}{c\phi^{\beta'}}.
\qed
\end{equation*}

\subsubsection{Proof of Lemma~\ref{lem.bal}}
\label{sec:prf.lem.bal}
Let $\gamma_{\phi}$ be the indecision level corresponding to classifier $Y^*$ with threshold $\tau^*_{\gamma}+\phi$, and let $\widehat{\gamma}_{\phi}$ be the indecision level corresponding to $\widehat{Y}$ with threshold $\tau^*_{\gamma}+\phi$. Define $\eta_{\mathrm{min}}=\eta\wedge(1-\eta)$ and $\eta_{\mathrm{max}}=\eta\vee(1-\eta)$, and note that
\begin{eqnarray*}
    R({Y^*}^{\tau^*_{\gamma}})-R({Y^*}^{\tau^*_{\gamma}+\phi}) &=& \tfrac1{1-\gamma}\mathbf{E}\eta_{\mathrm{min}} \mathbf{1}\{Y^*_{\gamma}\ne0\}
    -\tfrac1{1-\gamma_{\phi}}\mathbf{E}\eta_{\mathrm{min}} \mathbf{1}\{Y^*_{\gamma_{\phi}}\ne0\}\\
    &=& \tfrac1{1-\gamma}\mathbf{E}\eta_{\mathrm{min}} \big(\mathbf{1}\{Y^*_{\gamma}\ne0\}-\mathbf{1}\{Y^*_{\gamma_{\phi}}\ne0\}\big)
    +\big(\tfrac1{1-\gamma}-\tfrac1{1-\gamma_{\phi}}\big)\mathbf{E}\eta_{\mathrm{min}} \mathbf{1}\{Y^*_{\gamma_{\phi}}\ne0\}\\
    &=& \tfrac1{1-\gamma}\mathbf{E}\eta_{\mathrm{min}} \mathbf{1}\{\tau^*_{\gamma}<\eta_{\mathrm{max}}\le\tau^*_{\gamma}+\phi\}
    +\tfrac{(\gamma-\gamma_{\phi})}{(1-\gamma)(1-\gamma_{\phi})}\mathbf{E}\eta_{\mathrm{min}} \mathbf{1}\{Y^*_{\gamma_{\phi}}\ne0\}\\
    &\ge& \tfrac1{1-\gamma}\Big[\mathbf{E}\eta_{\mathrm{min}} \mathbf{1}\{\tau^*_{\gamma}<\eta_{\mathrm{max}}\le\tau^*_{\gamma}+\phi\}
    -(\gamma_{\phi}-\gamma)(1-\tau^*_{\gamma}-\phi)\Big]\\
    &=& \tfrac{(\gamma_{\phi}-\gamma)}{(1-\gamma)}\Big[\phi - \mathbf{E}\Big(\eta_{\mathrm{max}} - \tau^*_{\gamma} \,\big|\,\tau^*_{\gamma}<\eta_{\mathrm{max}} \le\tau^*_{\gamma}+\phi\Big)\Big].
\end{eqnarray*}
Condition (\ref{lbnds_eta_prob01})  implies that $\gamma_{\phi}-\gamma\gtrsim \phi^{\beta'}$ and $\phi-\mathbf{E}\big(\eta_{\mathrm{max}}-\tau^*_{\gamma} \,\big|\,\tau^*_{\gamma}<\eta_{\mathrm{max}}\le\tau^*_{\gamma}+\phi\big)\gtrsim \phi^{1+\beta' - \beta}$, where the second inequality follows from
\begin{align*}
\mathbf{E}\Big(\eta_{\mathrm{max}} - \tau^*_{\gamma} \,\big|\, \tau^*_{\gamma} < \eta_{\mathrm{max}} \le \tau^*_{\gamma} + \phi\Big) 
&\leq \frac{\frac{\phi}{2} \mathbf{P}(0 < \eta_{\mathrm{max}} - \tau^*_{\gamma} \le \phi/2) + \phi\left(\gamma_{\phi} - \gamma - \mathbf{P}(0 < \eta_{\mathrm{max}} - \tau^*_{\gamma} \le \phi/2)\right)}{\gamma_{\phi} - \gamma} \\
&\leq \phi - \frac{\phi\, \mathbf{P}(0 < \eta_{\mathrm{max}} - \tau^*_{\gamma} \le \phi/2)}{2\, \mathbf{P}(0 < \eta_{\mathrm{max}} - \tau^*_{\gamma} \le \phi)}.
\end{align*}
Consequently, $R({Y^*}^{\tau^*_{\gamma}})-R({Y^*}^{\tau^*_{\gamma}+\phi})\gtrsim\phi^{1+2\beta'-\beta}$, and we have derived the first bound of Lemma~\ref{lem.bal}.

Taking advantage of the fact that the threshold used by~${\widehat{Y}}^{\tau^*_{\gamma}+\phi}$ and~${Y^*}^{\tau^*_{\gamma}+\phi}$ is the same, and repeating the standard argument in \cite{Herbei_Wegkamp_06} while conditioning on~$\widehat{\eta}$, we derive that
\begin{eqnarray*}
R_{\widehat{\eta}}({\widehat{Y}}^{\tau^*_{\gamma}+\phi})-R({Y^*}^{\tau^*_{\gamma}+\phi})&\lesssim& \mathbf{P}\{ |\widehat{\eta}-\eta|>\epsilon\,|\,\widehat{\eta}\} +\epsilon\big[\mathbf{P}\big(|\tau^*_{\gamma}+\phi-\eta|\le\epsilon\big)+\mathbf{P}\big(|1-\tau^*_{\gamma}-\phi-\eta|\le\epsilon\big)\big].
\end{eqnarray*}
Thus, using $\epsilon\le\phi$ together with condition (\ref{lbnds_eta_prob01}) we arrive at the second bound of Lemma~\ref{lem.bal}. \qed

\subsection{Proof of Theorem \ref{thm:multi:lower_bound}}
     Let us consider a classifier $\tilde Y(X)$ such that $\mathbf{P}(\tilde Y = 0) = \gamma$ and let $A$ be the set where $\tilde{Y} = 0$. We have that
    \begin{align*}
    {\mathbf P}_Y(\tilde Y \neq Y | \tilde Y \neq 0 ) &= 1 - \frac{\sum_i p_i \mathbf{P}_i(\tilde Y = i) }{1 - \mathbf{P}(\tilde Y = 0)}\\
    & = 1 - \frac{\int_{A^c}\sum_i\mathbf{1}\{\tilde Y(x)=i\}p_if_i(x)\,d\mu(x)}{1-\gamma}.
    \end{align*}
For each $x$, the integrand is maximized by $\tilde Y = \underset{i}{\arg\max} (p_i f_i)$. Hence
$$
    {\mathbf P}_Y(\tilde Y \neq Y | \tilde Y \neq 0 ) \geq  1 - \frac{\int_{A^c}\underset{i}{\max}(p_if_i)\,d\mu}{1-\gamma}.
$$
Invoking Lemma \ref{lem:constrained:min} we get further that the above quantity is minimized for 
$$
A^* := \left\{x \; ; \; \underset{i}{\max} p_i f_i \leq \tau_\gamma \sum_i p_i f_i  \right\},
$$
such that $\mathbf{P}(A^*) = \gamma$ and $\tau_\gamma\in[1/K,1]$. The result and the expression for $Y^*$ follow.

\subsection{Proof of Theorem \ref{thm:asymp}}
The following bound for the tail of the Gaussian distribution will be useful for this proof. For all $t>0$, we have
$$
\frac{\exp(-t^2/2)}{\sqrt{2\pi}(t+1)}\leq\mathbf P(\xi\geq t)\leq\frac{\exp(-t^2/2)}{\sqrt{2\pi}\,t}.
$$
    We start with the case $1/2<c<1$:

    Remember that $\Delta = c\sqrt{2\log(1/\delta)}$. Let us choose $t = (1-c)\sqrt{2\log(1/\delta)}$. In that case
    $$
    \frac{\delta}{\sqrt{2\pi}(\sqrt{ 2\log(1/\delta)}+1) } \leq \mathbf{P}(\xi \geq \Delta + t) \leq \frac{\delta}{\sqrt{4\pi \log(1/\delta) }},
    $$
    and
    $$
    \frac{\delta^{(2c-1)^2}}{\sqrt{2\pi}((2c-1)\sqrt{ 2\log(1/\delta)}+1) } \leq \mathbf{P}(\xi \geq  \Delta - t ) \leq \frac{\delta^{(2c-1)^2}}{(2c-1)\sqrt{4\pi \log(1/\delta) }}.
    $$
    It follows that
    $$
    \frac{\mathbf{P}(\xi \geq \Delta + t)}{\mathbf{P}(\xi \geq   t- \Delta  )} \leq \frac{\delta}{\sqrt{4\pi \log(1/\delta) } \left(1 - \frac{\delta^{(2c-1)^2}}{\sqrt{2\pi}((2c-1)\sqrt{ 2\log(1/\delta)}) } \right)}.
    $$
     It is now clear that for small values of $\delta$ we have
     $$
\frac{\mathbf{P}(\xi \geq \Delta + t)}{\mathbf{P}(\xi \geq \Delta + t) + \mathbf{P}(\xi \geq t - \Delta )} \leq 
\delta.
     $$
     As a consequence
     \begin{equation}\label{eq:c>1/2:1}
              \gamma_\delta \leq \frac{\delta^{(2c-1)^2}}{(2c-1)\sqrt{4\pi \log(1/\delta) }} - \frac{\delta}{\sqrt{2\pi}(\sqrt{ 2\log(1/\delta)}+1) }.
     \end{equation}

     For $0<\varepsilon<1-c$, let us now choose $t = (1-c - \varepsilon)\sqrt{2\log(1/\delta)}$. In that case
    $$
    \frac{\delta^{(1-\varepsilon)^2}}{\sqrt{2\pi}((1-\varepsilon)\sqrt{2\log(1/\delta)}+1)}\leq\mathbf P(\xi\geq\Delta+t)\leq\frac{\delta^{(1-\varepsilon)^2}}{(1-\varepsilon)\sqrt{4\pi\log(1/\delta)}},
    $$
    and
    $$
    \frac{\delta^{(2c-1+\varepsilon)^2}}{\sqrt{2\pi}((2c-1+\varepsilon)\sqrt{ 2\log(1/\delta)}+1) } \leq \mathbf{P}(\xi \geq  \Delta - t ) \leq \frac{\delta^{(2c-1+\varepsilon)^2}}{(2c-1+\varepsilon)\sqrt{4\pi \log(1/\delta) }}.
    $$
    It follows that
    $$
    \frac{\mathbf P(\xi\geq\Delta+t)}{\mathbf P(\xi\geq t-\Delta)}\geq
    \frac{\delta^{(1-\varepsilon)^2}}
    {\sqrt{2\pi}((1-\varepsilon)\sqrt{2\log(1/\delta)}+1)
    \left[1-\frac{\delta^{(2c-1+\varepsilon)^2}}
    {\sqrt{2\pi}((2c-1+\varepsilon)\sqrt{2\log(1/\delta)}+1)}\right]}.
    $$
     It is now clear that for small values of $\delta$ we have
     $$
\frac{\mathbf{P}(\xi \geq \Delta + t)}{\mathbf{P}(\xi \geq \Delta + t) + \mathbf{P}(\xi \geq t - \Delta )} \geq 
\delta.
     $$
     As a consequence, for any $0<\varepsilon<1-c$, we get
          \begin{equation}\label{eq:c>1/2:2}
     \gamma_\delta \geq \frac{\delta^{(2c-1+\varepsilon)^2}}{\sqrt{2\pi}((2c-1+\varepsilon)\sqrt{2\log(1/\delta)}+1)}-\frac{\delta^{(1-\varepsilon)^2}}{(1-\varepsilon)\sqrt{4\pi\log(1/\delta)}}.
     \end{equation}
     Combining \eqref{eq:c>1/2:1} and \eqref{eq:c>1/2:2}, we conclude that if $\log(1/\gamma_\delta)\leq(2c-1)^2\log(1/\delta)$, then $\delta$-consistency is possible. On the other hand, if $\log(1/\gamma_\delta)\geq(1+\varepsilon)(2c-1)^2\log(1/\delta)$, then $\delta$-consistency is impossible.

      We will next deal with the case $0<c<1/2:$

    Remember that $\Delta = c\sqrt{2\log(1/\delta)}$. Let us choose $t = 1/(4c)\sqrt{2\log(1/\delta)}$. In that case
    $$
    \frac{\delta^{(c+1/(4c))^2}}{\sqrt{2\pi}((c+1/(4c))\sqrt{ 2\log(1/\delta)}+1) } \leq \mathbf{P}(\xi \geq \Delta + t) \leq \frac{\delta^{(c+1/(4c))^2}}{(c+1/(4c))\sqrt{4\pi \log(1/\delta) }},
    $$
    and
    $$
    \frac{\delta^{(c-1/(4c))^2}}{\sqrt{2\pi}((1/(4c)-c)\sqrt{2\log(1/\delta)}+1)}\leq\mathbf P(\xi\geq t-\Delta)\leq\frac{\delta^{(c-1/(4c))^2}}{(1/(4c)-c)\sqrt{4\pi\log(1/\delta)}}.
    $$
    It follows that
    $$
    \frac{\mathbf{P}(\xi \geq \Delta + t)}{\mathbf{P}(\xi \geq   t- \Delta  )} \leq \delta \frac{\sqrt{2\pi}((1/(4c)-c)\sqrt{ 2\log(1/\delta)}+1)}{(c+1/(4c))\sqrt{4\pi \log(1/\delta) }}.
    $$
     It is now clear that for small values of $\delta$ and any $c>0$ we have
     $$
\frac{\mathbf{P}(\xi \geq \Delta + t)}{\mathbf{P}(\xi \geq \Delta + t) + \mathbf{P}(\xi \geq t - \Delta )} \leq 
\delta.
     $$
     As a consequence
     \begin{equation}\label{eq:c<1/2:1}
              \gamma_\delta \leq 1 -  \frac{\delta^{(c-1/(4c))^2}}{\sqrt{2\pi}((1/(4c)-c)\sqrt{ 2\log(1/\delta)}+1) }   -   \frac{\delta^{(c+1/(4c))^2}}{\sqrt{2\pi}((c+1/(4c))\sqrt{ 2\log(1/\delta)}+1) }  .
     \end{equation}

     For a choice of $0<\varepsilon<1-4c^2$, let us now choose $t = (1-\varepsilon)/(4c)  \sqrt{2\log(1/\delta)}$, so that $t>\Delta$. In that case
    $$
    \frac{\delta^{(c+(1-\varepsilon)/(4c))^2}}{\sqrt{2\pi}((c+(1-\varepsilon)/(4c))\sqrt{ 2\log(1/\delta)}+1) } \leq \mathbf{P}(\xi \geq \Delta + t) \leq \frac{\delta^{(c+(1-\varepsilon)/(4c))^2}}{(c+(1-\varepsilon)/(4c))\sqrt{4\pi \log(1/\delta) }},
    $$
    and
    $$
    \frac{\delta^{(c-(1-\varepsilon)/(4c))^2}}{\sqrt{2\pi}(((1-\varepsilon)/(4c)-c)\sqrt{2\log(1/\delta)}+1)}\leq\mathbf P(\xi\geq t-\Delta)\leq\frac{\delta^{(c-(1-\varepsilon)/(4c))^2}}{((1-\varepsilon)/(4c)-c)\sqrt{4\pi\log(1/\delta)}}.
    $$
    It follows that
    $$
    \frac{\mathbf P(\xi\geq\Delta+t)}{\mathbf P(\xi\geq t-\Delta)}\geq
    \delta^{1-\varepsilon}\,
    \frac{((1-\varepsilon)/(4c)-c)\sqrt{2\log(1/\delta)}}
    {(c+(1-\varepsilon)/(4c))\sqrt{2\log(1/\delta)}+1}.
    $$
     It is now clear that for small values of $\delta$ we have
     $$
\frac{\mathbf{P}(\xi \geq \Delta + t)}{\mathbf{P}(\xi \geq \Delta + t) + \mathbf{P}(\xi \geq t - \Delta )} \geq 
\delta.
     $$
     As a consequence, for small values of $\varepsilon>0$, we get
          \begin{equation}\label{eq:c<1/2:2}
     \gamma_\delta \geq 1 - \frac{\delta^{(c-(1-\varepsilon)/(4c))^2}}{\sqrt{2\pi}(((1-\varepsilon)/(4c)-c)\sqrt{ 2\log(1/\delta)}) }   -   \frac{\delta^{(c+(1-\varepsilon)/(4c))^2}}{\sqrt{2\pi}((c+(1-\varepsilon)/(4c))\sqrt{ 2\log(1/\delta)}) }  .
     \end{equation}
     Combining \eqref{eq:c<1/2:1} and \eqref{eq:c<1/2:2}, we conclude that if $\log(1/(1-\gamma_\delta))\geq(c-1/(4c))^2\log(1/\delta)$, then $\delta$-consistency is possible. On the other hand, if $\log(1/(1-\gamma_\delta))\leq(1-\varepsilon)(c-1/(4c))^2\log(1/\delta)$, then $\delta$-consistency is impossible.

\end{document}